\documentclass{article}

\usepackage{arxiv}

\usepackage{booktabs}
\usepackage{lineno,hyperref}
\hypersetup{colorlinks, citecolor=blue, linkcolor=red, urlcolor=green}
\usepackage[flushleft]{threeparttable}
\usepackage{compactbib}
\usepackage{amsmath,amssymb,amsfonts,bm}
\usepackage{arydshln}
\usepackage{graphics}
\usepackage{graphicx}
\usepackage{multirow}
\usepackage{caption}
\usepackage{subcaption}
\usepackage{subfloat}
\usepackage{acronym}
\usepackage{adjustbox}
\usepackage[sort&compress,numbers]{natbib}

\usepackage{algorithm}
\PassOptionsToPackage{noend}{algpseudocode}
\usepackage{algpseudocode}

\algdef{SE}[DOWHILE]{Do}{doWhile}{\algorithmicdo}[1]{\algorithmicwhile\ #1}%

\errorcontextlines\maxdimen
\usepackage{url}

\title{Data-driven discovery of interpretable Lagrangian of stochastically excited dynamical systems}

\date{} 					

\author{Tapas Tripura\\
	Department of Applied Mechanics,\\
    Indian Institute of Technology Delhi, \\
    Hauz Khas, Delhi 110016, India\\
	\texttt{tapas.t@am.iitd.ac.in} \\
	\And
	Satyam Panda\\
	Department of Mechanical Engineering\\
    Indian Institute of Technology Guwahati \\
    Assam 781039, India\\
	\texttt{panda18@iitg.ac.in} \\
    \And
 	Budhaditya Hazra \\
	Department of Civil Engineering,\\
    Indian Institute of Technology Guwahati,\\
    Assam 781039, India\\
	\texttt{budhaditya.hazra@iitg.ac.in} \\
    \And
	Souvik Chakraborty \\
	Department of Applied Mechanics,\\
    Indian Institute of Technology Delhi, \\
    Hauz Khas, Delhi 110016, India\\
	\texttt{souvik@am.iitd.ac.in} \\
}

\renewcommand{\headeright}{Tripura et al.}
\renewcommand{\undertitle}{}
\renewcommand{\shorttitle}{Discovery of Lagrangian of Stochastic Dynamical Systems}

\begin{document}
\maketitle

\begin{abstract}
Exploring the intersection of deterministic and stochastic dynamics, this paper delves into Lagrangian discovery for conservative and non-conservative systems under stochastic excitation. Traditional Lagrangian frameworks, adept at capturing deterministic behavior, are extended to incorporate stochastic excitation. The study critically evaluates recent computational methodologies for learning Lagrangians from observed data, highlighting the limitations in interpretability and the exclusion of stochastic excitation. To address these gaps, an automated data-driven framework is proposed for the simultaneous yet uncoupled discovery of Lagrange densities and the volatility function of stochastic excitation by leveraging the sparse regression. This novel framework offers several advantages over existing approaches. Firstly, it provides an interpretable description of the underlying Lagrange density, allowing for a deeper understanding of system dynamics under stochastic excitations. Secondly, it identifies the interpretable form of the generalized stochastic force, addressing the limitations of existing deterministic approaches. Additionally, the framework demonstrates robustness and versatility through numerical case studies encompassing both stochastic differential equations (SDEs) and stochastic partial differential equations (SPDEs), with results showing almost exact approximations to true system behavior and minimal relative error in derived equations of motion.
\end{abstract}

\keywords{Lagrangian discovery\and  Stochastic dynamics\and  Data-driven framework \and  Sparse regression\and  Interpretability}

\section{Introduction}
In the domain of modeling dynamical systems, the foundational frameworks of Lagrangian and Hamiltonian mechanics have proven vital for revealing the behavior of a broad spectrum of physical phenomena \cite{mann2018lagrangian,materassi2020stochastic,finzi2020simplifying}. These deterministic formulations traditionally excel in capturing the regular evolution of systems governed by conservation laws and are also adept at describing the dynamics of non-conservative systems. Expanding on this concept, the Lagrangian's versatility lies in its ability to streamline the representation of a system's dynamics into a single mathematical expression \cite{calkin1996lagrangian}. In contrast to equation discovery approaches, which elicit functional dependencies among system variables, the use of Lagrangian mechanics allows researchers to uncover fundamental governing principles, such as the conservation of energy and momentum, directly from the system's dynamics. This distinction highlights the unique ability of Lagrangian-based methods to provide deeper insights into the underlying physics governing system behavior. However, the complexity of real-world dynamics often extends beyond the purview of deterministic models, necessitating an exploration into the stochastic domain \cite{materassi2020stochastic}. Many natural processes exhibit inherent randomness, challenging our conventional understanding and prompting the incorporation of stochastic excitation. The Lagrangian, traditionally a deterministic concept, can be extended to incorporate stochastic excitation, allowing for a more realistic representation of systems subject to random influences \cite{eyink2020stochastic}. Stochastic Euler-Lagrange differential equations \cite{bou2009stochastic} and stochastic Hamiltonian differential equations \cite{panda2023geometry,holm2018stochastic} serve as a potent tool for modeling such phenomena, encapsulating both the deterministic dynamics dictated by the Lagrange density and the stochastic perturbations contributing to the intrinsic uncertainty in a system's responses. The application of the stochastic Lagrangian model extends beyond classical mechanics and finds relevance in the domain of quantum field theory \cite{lindgren2019quantum} and astrophysics \cite{trac2003primer,li2023lagrangian}. In contrast to conventional equation discovery, which focuses on identifying explicit mathematical relationships within data, the pursuit of the Lagrange density delves deeper into uncovering the underlying principles governing the dynamics of stochastic Hamiltonian systems. By elucidating the Lagrange density from observed data, this research aims to reveal the fundamental governing structure that encapsulates both deterministic dynamics and stochastic influences, providing insights into system behavior.

With the ongoing progress in data-driven methodologies, researchers are actively exploring approaches to directly learn the Lagrangian of systems from observational data. A notable development in this domain is the Hamiltonian neural network (HNN) \cite{toth2019hamiltonian,greydanus2019hamiltonian} and Hamiltonian graph neural networks \cite{sanchez2019hamiltonian}, which leverages the Hamiltonian principles to infer the Hamiltonian directly from data. However, a limitation arises as the Hamiltonian formalism mandates that the system's coordinates must be canonical, posing restrictions that may be impractical for real-world systems. In response to this challenge, deep Lagrangian networks (DeLaN) \cite{lutter2019deep} and Lagrangian neural networks (LNN) \cite{cranmer2020lagrangian} have emerged as a viable alternative. Additionally, efforts have been directed towards overcoming constraints, leading to the development of Constrained Hamiltonian Neural Networks (CHNNs) and Constrained Lagrangian Neural Networks (CLNNs), which enable the direct learning of Hamiltonians and Lagrangians in Cartesian coordinates \cite{gruver2022deconstructing}. Another noteworthy approach is the Lagrangian graph neural network (LGNN), designed specifically for learning the Lagrangian in the context of rigid body dynamics \cite{bhattoo2022learning}.

Recent advancements in Lagrangian discovery from observed data signify a dynamic evolution in computational methodologies. Nevertheless, a noticeable limitation emerges from the lack of interpretability within the identified Lagrangian models, presenting a substantial hindrance to their broad adoption and adaptability. The primary challenge stems from the inability to deduce explicit governing equations directly from these discovered Lagrangians, thereby constraining our capacity to attain profound insights into the fundamental dynamics of the systems under consideration. 
An alternative route involves utilizing recently developed library-based methodologies to discover Lagrangian from data \cite{tripura2023bayesian,tripura2024discovering}. The core idea here is to select suitable candidates from an extensive library of possible basis functions. The primary advantage of library-based approaches, originally proposed for discovering ordinary \cite{schmidt2009distilling,brunton2016discovering,nayek2021spike,nayek2022equation,wentz2023derivative}, partial \cite{schaeffer2017learning,rudy2017data,chen2021physics,flaschel2021unsupervised,joshi2022bayesian,more2023bayesian}, and stochastic differential equations from data \cite{boninsegna2018sparse,li2021data,tripura2023sparse,tripura2023probabilistic,mathpati2024discovering}, resides in its interpretability and perpetual generalization. Naturally, employing a library-based approach for Lagrangian discovery also offers interpretability and can facilitate the discovery of conservation law and governing equations through the application of the Noether theorem \cite{Noether1}.
However, a notable drawback of these approaches is their exclusion of stochastic excitation, potentially undermining the accuracy and applicability of the discovered models by failing to account for the inherent randomness and uncertainty present in real-world systems. Without incorporating stochastic influences, the discovered models may provide incomplete representations of system dynamics, leading to inaccuracies in predictions and suboptimal performance in control and optimization strategies. This omission becomes particularly critical when attempting to discover the Lagrange density of a stochastically excited dynamical system. Addressing this limitation necessitates future advancements in Lagrangian discovery frameworks that explicitly account for stochastic influences, ensuring a more accurate and comprehensive representation of the dynamics inherent in stochastically excited systems.

This paper addresses the identified limitations by introducing an automated framework for the simultaneous but uncoupled discovery of the Lagrange density and the volatility function of stochastic excitation from observational data. First and second-order statistical moments are employed to obtain an expression for Lagrangian density and volatility function in terms of measurement data.
Overall, this work aims to address several key challenges; firstly, it focuses on developing an automated framework for the discovery of interpretable Lagrangian density for stochastically excited systems. Secondly, it automates the discovery of volatility functions in stochastic systems, a natural aspect that is neglected in the deterministic model discovery frameworks. Thirdly, leveraging the stochastic Euler-Lagrangian formalism on the discovered Lagrangian density and volatility function, it uncovers highly interpretable SDE or SPDE descriptions of underlying stochastic systems. 
To enhance interpretability and directly extract the essential dynamics from data, we employ the sequential threshold least-squares \cite{kaheman2020sindy}. Further, the automated framework is designed to operate efficiently with the requirement of the only ensemble of the system output responses. This streamlined process enhances the practicality, adaptability, and accessibility of the proposed framework to real-world scenarios where obtaining input force information is intractable.
The novel contributions of the proposed framework can be encapsulated in the following points: (I) This method systematically uncovers the governing equations while maintaining interpretability, (II) This approach provides a more comprehensive understanding of the system's behavior under random excitations, and (III) The proposed framework is generalized for systems with states featuring both temporal and/or spatial derivatives, i.e., Stochastic Differential Equations (SDEs) and Stochastic Partial Differential Equations (SPDEs), expanding its applicability to a broader range of complex dynamical systems. 
The efficacy of the proposed framework is investigated on six stochastic dynamical systems (involving both SDEs and SPDEs). The rediscovered equations closely adhere to the true SDEs/SPDEs with approximately 1\% relative error in most cases, indicating promising applications to real-world problems where stochastic fluctuations are a dominating factor. 

The rest of the paper is structured as follows: Section \ref{sec:problem} outlines the problem formulation, while Section \ref{sec:lagrangian_discovery} introduces the proposed automated data-driven framework for discovering the Lagrange density and the volatility of stochastic excitation from observational data. In Section \ref{sec:numerical}, numerical experiments are conducted to highlight the applicability of the proposed data-driven framework. Section \ref{sec:conclusion} concludes this study by highlighting key features of the proposed framework.

\section{Problem Formulation}\label{sec:problem}
In the investigation of continuous dynamical systems under stochastic excitation, the sensors are placed at finite locations to collect the data. It is assumed that these sensors proficiently capture the spatial variations, each meticulously recording the temporal data. This conceptualization permits the representation of the system's Lagrangian based on the sensor measurements. Consequently, the continuum system is modeled as an $n$-particle system, wherein each particle is analogous to the position of a sensor. These complex systems respond to stochastic influences, and the task is to extract valuable information about their underlying dynamics. In practical scenarios, we are provided with an ensemble of the responses of these systems ${\bm{U}}(t) = \left[u_i(t), (u_x (t))_i, (u_t(t))_i \right]_{j=1}^{N} \in \mathbb{R}^{N_t \times n \times N}$ with $u_t$ and $u_x$ being temporal and spatial derivative of $u$ and $i = 1,\ldots,n$ represents the number of particles and $N$ represents the ensemble count. In the upcoming discussions, we'll simplify our notation by dropping the explicit time dependence $t$ and represent system states as $\left(u, u_x, u_t\right)$. To this end, we define the Lagrangian as follows:
\begin{equation}
    \mathcal{L}\left(u, u_t, u_x\right) = T\left(u,u_t\right) - V\left(u,u_x\right),
\end{equation}
where $T\left(u,u_t\right) \in \mathbb{R}$ represents the kinetic energy and $V\left(u,u_x\right) \in \mathbb{R}$ represents the potential energy.
The action integral, denoted as $\mathcal{A}$, involves integrating the Lagrangian density along the system's motion throughout a specified time duration and a spatial region. For a discrete spatial domain with $n$ number of particles, considering Lagrangian at $i^{th}$ particle as $\mathcal{L}_i$ at $x = l$ and considering the fixed interval $t \in \left[t_o,t_f\right]$, we define the action integral as follows,
\begin{equation} \label{act}
    \begin{aligned}
        \mathcal{A} &= \int_{t_0}^{t_f} \left( \mathcal{L}_i - \bar{\mathcal{S}}_i \right) \, dt,
    \end{aligned}
\end{equation}
where $\bar{\mathcal{S}}_i$ is the external work done for $i^{th}$ particle, which is for a Wiener noise $\mathcal{W}(l,t) \sim \mathcal{N}(0,t)$ is defined as $ \sigma_{i}\left(u\left(l,t\right) \right) \circ \dot{\mathcal{W}_t}$. Here, $\dot{\mathcal{W}}_t \sim \mathcal{N}(0,1)$ is the white noise, defined as the derivative of Wiener noise $\mathcal{W}(l,t)$ with respect to time $t$, and $\circ$ is the Stratonovich operator, which ensures consistent treatment of noise terms in calculus operations. The preference for Stratonovich or Itô is subjective, emphasizing the need for context-specific considerations in selecting the appropriate mathematical framework. For example, in It\^{o} calculus \cite{ito1944109}, conventional integration by parts and the chain rule (Newton–Leibniz rule) are inapplicable. As a consequence, in fields like mathematical finance \cite{tankov2003financial} or population biology \cite{turelli1977random}, stochastic Itô integrals integral are used. As an alternative, Stratonovich integral \cite{chechkin2014marcus} is often preferred in classical mechanics due to its consistency with deterministic calculus rules and its ability to preserve symmetries and physical structure. It is also valuable in scenarios dealing with a Wong–Zakai small correlation time, especially when noise arises as a continuous approximation of a rapidly fluctuating force. Using the definition of external work done, the Eq. \eqref{act} can be rephrased as,
\begin{equation}
    \begin{aligned}
        \mathcal{A} &= \int_{t} \left( \mathcal{L}_i -  \left(\sigma_{i}\left(u\left(l,t\right) \right) \circ \dot{\mathcal{W}_t}\right) \right) \, dt .
    \end{aligned}
\end{equation}
For the curve $u(l,t):\left[t_0, t_f\right] \mapsto \mathbb{R}^n$, the variation $u(l,t)$ is characterized by an $\epsilon$-parametrized family of curves $u^\epsilon(l,t)$ within $\mathbb{R}^n$, where $\epsilon \in (-c, c)$ for some positive constant $c$, and $u^0(l,t)=u(l,t)$. As the end points are fixed, $u^\epsilon\left(l,t_0\right)=u\left(l,t_0\right)$ and $u^\epsilon\left(l,t_f\right)=u\left(l,t_f\right)$. Expanding the $u^\epsilon(l,t)$ and taking an infinitesimal variation, we get
\begin{equation}
    \left. \frac{d}{d\epsilon} u^\epsilon \right|_{\epsilon = 0} = \delta u .
\end{equation}
The same process can be carried out for obtaining the infinitesimal variation of the curve $u_t(l,t)$ as,
\begin{equation}
    \begin{aligned}
        \left. \frac{d}{d\epsilon} u_t^\epsilon \right|_{\epsilon = 0} = \delta u_t .
    \end{aligned}
\end{equation}
The principle of stationary action states that $\delta \mathcal{A} = 0$, where $\delta \mathcal{A}$ denotes its first variation with respect to $u$, $u_t$, and $u_x$.
The variation of the action is expressed as,
\begin{equation}\label{eq:variation}
    \begin{aligned}
        \delta \mathcal{A} &= \int_{t} \left( \frac{\partial \mathcal{L}_i}{\partial u} \delta u + \frac{\partial \mathcal{L}_i}{\partial u_t} \delta u_t \right) \, dt - \int_{t} \left(\frac{\partial \sigma_{i}}{\partial u} \circ \dot{\mathcal{W}_t} \right) \delta u \, dt , \\
        &= \int_{t} \left( \frac{\partial \mathcal{L}_i}{\partial u} - \frac{\partial}{\partial t} \frac{\partial \mathcal{L}_i}{\partial u_t} \right) \delta u \, dt - \int_{t} \left(\frac{\partial \sigma_{i}}{\partial u} \circ \dot{\mathcal{W}_t} \right) \delta u \, dt .
    \end{aligned}
\end{equation}
As white noise is nowhere differentiable, it is often represented using the increments of the Wiener process. Despite its lack of differentiability in the classical sense, white noise is associated with the concept of mean square derivative of Wiener process, where the mean square derivative provides insight into its behavior by considering the average rate of change over time through the notion of quadratic variation. Using the properties of mean squared derivative, the differential Wiener process is defined as $d\mathcal{W} \sim \mathcal{N}(0, dt)$, i.e., Gaussian white noise (GWN) process with variance $dt$. Incorporating this result along with the natural boundary conditions in Eq. \eqref{eq:variation}, we obtain the stochastic Euler-Lagrange equation:
\begin{equation}\label{SEL}
    \frac{\partial \mathcal{L}_i}{\partial u} - \frac{\partial}{\partial t} \frac{\partial \mathcal{L}_i}{\partial u_t} = \frac{\partial \sigma_{i}}{\partial u} \circ \dot{\mathcal{W}_t}
\end{equation}
Eq. (\ref{SEL}) represents Stratanovich SDE and is transformed into Ito counterpart aided by Wong-Zakai \cite{williams2012stochastic} correction terms as follows,
\begin{equation}\label{SEL1}
    \frac{\partial \mathcal{L}_i}{\partial u} - \frac{\partial}{\partial t} \frac{\partial \mathcal{L}_i}{\partial u_t} - \frac{1}{2}\frac{\partial \sigma_{i}}{\partial u} \frac{\partial^2 \sigma_{i}}{\partial u^2}  = \frac{\partial \sigma_{i}}{\partial u} \dot{\mathcal{W}_t}
\end{equation}
For the systems with additive noise, Eq. (\ref{SEL1}) is reduced to the following form,
\begin{equation}\label{SEL2}
    \frac{\partial \mathcal{L}_i}{\partial u} - \frac{\partial}{\partial t} \frac{\partial \mathcal{L}_i}{\partial u_t} = \frac{\partial \sigma_{i}}{\partial u} \dot{\mathcal{W}_t}
\end{equation}
The aim of this work is to learn an interpretable description of the Lagrangian of $n$-particle dynamical systems, denoted by $\mathcal{L}_{i} \in \mathbb{R}$ from sensor measurements. It also aims to discover the volatility function ${\partial \sigma_{i}}/{\partial u}$ from the observed responses. 
The discovered Lagrangian $\mathcal{L}_i$ serves as a key mathematical backbone that encapsulates the system's kinetic and potential energy contributions and the volatility function ${\partial \sigma_{i}}/{\partial u}$ represents the influence of stochastic force on the system dynamics.
For automated discovery of the interpretable expressions, the proposed framework leverages sparse lest-squared regression, which is discussed in brief in the following sections. We again highlight that the proposed requires only an ensemble of system responses, which can be obtained by assuming the underlying physical system to be ergodic.

\section{Lagrangian discovery}\label{sec:lagrangian_discovery}
This section introduces an automated framework designed for the systematic discovery of the Lagrangian density and diffusion, leveraging the observed responses of the dynamical system. The approach here is grounded in the principles of inverse problem-solving, wherein the aim is to infer the underlying dynamics of the system from its observable behavior. Towards learning the Lagrange density, the assumption here is that the Lagrangian can be represented as a weighted sum of specific candidate basis functions as,
\begin{equation}\label{eq:lagrange}
    \begin{aligned}
        \mathcal{L}_{j} &= \sum_{j=1}^n \sum_{i=1}^m c_{ij} D_i\left(u, u_t, u_x\right), \;\; j = 1,\ldots,n,
    \end{aligned}
\end{equation}
where ${\bm{C}}_i = \left\{c_1,\ldots, c_m \right\}$ represents the coefficients (the system parameters) of candidate basis functions $\bm{D} = \{{D}_1, \ldots, {D}_m \}$ for $i^{th}$ particle. Here, $m$ denotes the number of candidate basis functions. These candidate basis functions are chosen from a symbolic library encompassing various energy potentials. This symbolic library includes functions like polynomials and harmonics applied to the system's states and polynomials of their differences. As part of this methodology, we have also incorporated the concept of stochastic discovery, demonstrating how these dictionaries are practically applied to systems described by Stochastic Differential Equations (SDEs) and Stochastic Partial Differential Equations (SPDEs). This showcases the adaptability and effectiveness of our framework, particularly in handling stochastically excited dynamical systems.
For further understanding, consider the dictionaries employed in this work for the stochastically excited systems as,
\begin{equation}
    \begin{aligned}
        \boldsymbol{D}\left(u, u_t, u_x\right)= & {\left[\begin{array}{llllllll}
        1 & \mathrm{P}^{\rho_1}(u) & \mathrm{P}^{\rho_2}(u_t) & \mathrm{P}^{\rho_3}(u_x) & \operatorname{Trig.}(u) & \operatorname{Trig.}(u_t) & \operatorname{Trig.}(u_x) & \ldots \\
        \end{array}\right], }
    \end{aligned}
\end{equation}
Here, $\mathrm{P}^{*}(\cdot)$ denotes the polynomial function with $\rho_1, \rho_2, \rho_3 \in \mathbb{R}$ as the degree of the polynomial, Trig.($\cdot$) denotes trigonometric functions like Sine and Cosine, and $u_t, u_x$ denote partial derivatives with respect to time $t$ and space $x$. The symbolic library $\boldsymbol{D}\left(u, u_t, u_x\right)$ is evaluated on the state measurements, where each column of $\boldsymbol{D}\left(u, u_t, u_x\right)$ represents a possible candidate to be included in the final model of the Lagrangian. The type and number of total basis functions need to be carefully chosen. In a compact matrix notation, the linear regression problem in Eq. \eqref{eq:lagrange} can be represented as,
\begin{equation}\label{disc}
    \mathcal{L}_{i}\left(u, u_t, u_x\right) = \sum_{i=1}^{n} {\bm{D}}\left(u, u_t, u_x\right)\bm{C}_{i}^{\top},  \;\; i = 1,\ldots,n.
\end{equation}
Once the Lagrangian $\mathcal{L}_{i}$ at $n$-particles are discovered, the complete Lagrangian of a $n$-particle system can be constructed as the linear sum of the Lagrangian of individual particles. 
In the context of stochastic Lagrangian systems, the inherent unpredictability and randomness of the dynamic processes pose a significant challenge to understanding their behavior. To address this, statistical measures become crucial for comprehensively characterizing the system's behavior. Moments, such as mean and variances, offer insights into the distributional properties of the system's trajectories within the ensemble. To that end, the subsequent subsections delve into the discovery of Lagrangian functional and volatility of the stochastic excitation.

\subsection{Discovery of Lagrangian functional}\label{sub1}
Leveraging the theory of statistical measures, the expectation of the stochastic Euler-Lagrange equation from Eq. (\ref{SEL}) can be written as, 
\begin{equation} \label{exp_SL}
    \mathbb{E}\left[ \frac{\partial \mathcal{L}_i}{\partial u} - \frac{\partial}{\partial t} \frac{\partial \mathcal{L}_i}{\partial u_t}\right] = \mathbb{E}\left[\frac{\partial \sigma_{i}}{\partial u} \frac{d \mathcal{W}_t}{dt} \right] , \;\; j = 1,\ldots,n,
\end{equation}
where $\mathcal{L}_i$ represents Lagrangian at $i^{th}$ particle at $x = l$ and $\mathbb{E}[\cdot]$ denotes the expectation operation. Let \((\Omega, \mathcal{F}, (\mathcal{F}_t)_{t \geq 0}, P)\) be a stochastic basis with a filtration \((\mathcal{F}_t)_{t \geq 0}\). 
A standard Wiener process $\mathcal{W}_t$ adapted to this filtration has normally distributed increments $\Delta \mathcal{W}_t \sim \mathcal{N}(0,\Delta t)$. The white noise $d\mathcal{W}_t$ is defined as the limit in the distribution of a sequence $\Delta \mathcal{W}_t^m$ such that $d\mathcal{W}_t = \lim_{m \to \infty} \Delta \mathcal{W}_t^m$. 
This definition ensures that $d\mathcal{W}_t$ is adapted to the filtration \((\mathcal{F}_t)_{t \geq 0}\), aligning with the constant intensity and uncorrelated increment properties of white noise. Since $\Delta \mathcal{W}_t$ has zero mean, the term $\mathbb{E}\left[({\partial \sigma_{i}}/{\partial u}) ({d \mathcal{W}_t}/{dt}) \right]$ on the right-hand side of the above equation vanishes. Therefore, the Eq. \eqref{exp_SL} simplifies to,
\begin{equation} \label{exp_SL1}
    \left(\frac{\partial }{\partial u} - \frac{\partial}{\partial t} \frac{\partial }{\partial u_t} \right) \mathbb{E}\left[\mathcal{L}_i\right] = 0 , \;\; j = 1,\ldots,n.
\end{equation}
Then utilizing Eq. (\ref{disc}), the above equation can be rewritten for the $i^{th}$ particle as,
\begin{equation}\label{eq:regression_zero}
    \left(\frac{\partial }{\partial u} - \frac{\partial}{\partial t} \frac{\partial }{\partial u_t} \right) \mathbb{E}\left[{\bm{D}}\right]{\bm{C}}_{i}^{\top} = 0 , \;\; j = 1,\ldots,n.
\end{equation}
With this setup, our aim is to identify an optimal set of basis functions from the library ${\bm{D}}$ and determine its corresponding coefficients ${\bm{C}}_i$ that provide the most precise representation of the Lagrangian functional $\mathcal{L}_i$ from the measured data. This makes the problem a variable selection cum parameter estimation problem.
However, we notice that the solution to Eq. \eqref{eq:regression_zero} has a trivial solution. Therefore, before proceeding to the variable selection problem, we perform a non-trivial treatment of the Eq. \eqref{eq:regression_zero}.
In particular, we extract the kinetic energy basis from the library and denote it as ${D}_{u_t^2}$. We use the extracted function ${D}_{u_t^2}$ as the label vector for sparse regression. With the extraction, the dimension of the basis functions in the actual feature directory ${\bm{D}}$ becomes $m-1$. We denote the reduced dictionary as ${\bm{D}}_{-u_t^2} \in \mathbb{R}^{m-1}$. This leads to the following sparse regression equation,
\begin{equation} \label{exp_SL2}
    \left(\frac{\partial }{\partial u} - \frac{\partial}{\partial t} \frac{\partial }{\partial u_t} \right) \mathbb{E}\left[{\bm{D}}_{u_t^2}\right] = \left(\frac{\partial }{\partial u} - \frac{\partial}{\partial t} \frac{\partial }{\partial u_t} \right) \mathbb{E}\left[{\bm{D}}_{-u_t^2}\right]{\bm{C}}_i^{\top} + \mathcal{E} , \;\; j = 1,\ldots,n,
\end{equation}
where ${\bm{D}}_{-u_t^2}$ is the candidate Lagrangian library without the column of kinetic energy basis and $\mathcal{E} \in \mathbb{R}$ is the model mismatch error. The precision of the identified Lagrangian is entirely contingent upon the selection and adaptability of the library functions. The considered library has the capacity to encompass a broad spectrum of energy potentials and kinetic energy functions dependent on both states. In Eq. \eqref{exp_SL2}, the only unknown parameter is the coefficient vector of the Lagrangian library. 
Amidst the absence of prior information, the library may encompass numerous basis functions. However, according to prior knowledge \cite{tripura2024discovering}, the general Lagrangian is comprised of only a limited subset of these functions, where the majority of entries in the vector ${\bm{C}}_i$ are expected to exhibit a zero value. 
An iterative process is then executed on the remaining non-zero entries of ${\bm{C}}_i$, with the cardinality of non-zero coefficients diminishing with each iteration. The termination criterion is set to the convergence of non-zero coefficients. Towards identifying the sparse solution vector ${\bm{C}}_i$, we leverage the sequential threshold least-square regression technique, as proposed in \cite{tripura2024discovering}. Utilizing repeated measurements to account for the intrinsic stochastic nature of the excitation process, our objective is to determine the least-squares solution for the coefficients ${\bm{C}}_{i}$ and subsequently applying a sparsity constant $\lambda$ to threshold the solution vector. This procedure is iteratively repeated until the non-zero coefficients converge. In each iteration, achieving sparsity involves imposing a penalty on the solution through its $\text{L}_1$-norm in the following manner:
\begin{equation} \label{STLS}
    \begin{array}{ll}
         {\bm{C}}_{i} = & \arg. \mathop {\min } \limits_{{\bm{C}}^\star_{i}} \left\| \left(\frac{\partial }{\partial u} - \frac{\partial}{\partial t} \frac{\partial }{\partial u_t} \right) \mathbb{E}\left[{\bm{D}}_{-u_t^2}\right]{\bm{C}}_i^{\star \top} - \left(\frac{\partial }{\partial u} - \frac{\partial}{\partial t} \frac{\partial }{\partial u_t} \right) \mathbb{E}\left[{\bm{D}}_{u_t^2}\right]\right\|_2  
          + \lambda \left\|{\bm{C}}^\star_i \right\|_1, \, \forall \, i - 1,\ldots,n 
    \end{array}
\end{equation}

where, $\| \cdot \|_1$ and $\| \cdot \|_2$ represents L$_1$ and L$_2$ norms, respectively. Sequential Threshold Least Square Regression employs $n$-distinct search to obtain a sparse matrix $\Theta = \{ {\bm{C}}_1^{\top}, {\bm{C}}_2^{\top}, \ldots, {\bm{C}}_n^{\top} \}$, each representing a sparse vector associated with a specific particle. 
Once the sparse matrix is obtained, the Lagrange density of each particle will contain the sparsely identified Lagrangian basis functions and the extracted kinetic energy basis. The discovered parametric model is finally obtained from the identified basis functions as follows,
\begin{equation} \label{Lagg}
    \mathcal{L}_i = {\bm{D}}_{-u_t^2}{\bm{C}}_{i}^{\top} + {\bm{D}}_{u_t^2}, \;\; i =1,\ldots,n.
\end{equation}
Finally, with the Lagrangian densities of $n$ particles, the final Lagrangian functional of complete systems is obtained as a linear sum of the Lagrangian density of individual particles as $\mathcal{L}_{s} = \sum_{i=1}^n \mathcal{L}_{i}$. 
Algorithm \ref{algo1} provides the implementation of the proposed framework for the discovery of the Lagrangian functional. 
Upon discovering the Lagrangian density, the subsequent step involves determining the diffusion of the stochastic excitation.  
\begin{algorithm}[H]
	\small
    \caption{Lagrangian discovery for stochastically excited system}
	\begin{algorithmic}
    \Require{Observed responses: $\{u_i(t), (u_x(t))_{i}, (u_t(t))_{i}\}^{N}_{j=1} \in \mathbb{R}^{N_t \times n \times N}$ and sparsification constant $\lambda \in \mathbb{R}$.} 
        \State {Select $m$ candidate energy potential functions.}
        \State {Construct a symbolic library $\bm{D}\left( u, u_x, u_t \right) \in \mathbb{R}^{m}$.}
        \For{$i = 1$ to $n$}\Comment{Iteration over particles}
        \State {Compute the derivatives $\frac{\partial \bm{D}}{\partial u}$ and $\frac{\partial \bm{D}}{\partial u_{t}}$.} 
        \State {Evaluate the derivatives over observations $\{ \frac{\partial \mathbf{D}}{\partial u}, \frac{\partial \mathbf{D}}{\partial u_{t}} \}_{j=1}^{N}$.} 
        \State {Compute the time derivative of momentum $\{ \frac{d}{dt} \bigl( \frac{\partial \mathbf{D}}{\partial u_{t}} \bigr) \}_{j=1}^{N}$.} \Comment{Refer Eq. (\ref{exp_SL1})} 
        \State {Take expectation and compute the Euler Lagrange library: $\frac{d}{dt} \bigl( \frac{\partial \mathbb{E}[\mathbf{D}]}{\partial u_{t}} \bigr) - \frac{\partial \mathbb{E}[\mathbf{D}]}{\partial u}$.} \Comment{Refer Eq. (\ref{exp_SL1})} 
        \State{Extract $\mathbb{E}[{\bm{D}}_{u_t^2}]$ and $\mathbb{E}[{\bm{D}}_{-u_t^2}]$.}
        \State {Construct the regression equation.} \Comment{Refer Eq. (\ref{exp_SL2})} 
        \State {Do sequentially threshold least-squares to obtain the coefficients ${\bm{C}}_i$.} \Comment{Refer Algorithm \ref{algo}} 
        \State Obtain the Lagrangian of $i^{\text{th}}$ particle: $\mathcal{L}_i = \mathbb{E}[{\bm{D}}_{-u_t^2}]{\bm{C}}_{i}^{\top} + \mathbb{E}[{\bm{D}}_{u_t^2}].$ \Comment{Refer Eq. (\ref{Lagg})}
        \EndFor
    \Ensure{Lagrangian of complete system: $\mathcal{L}_{s}=\sum_{i=1}^{n} \mathcal{L}_{i}$.}
	\end{algorithmic}
	\label{algo1}
\end{algorithm}
\subsection{Discovery of volatility of the stochastic excitation}
In this context, the utilization of the second-order moment proves instrumental. Taking the second order moment of Eq. (\ref{SEL}), one can get,
\begin{equation} \label{var_SL}
    {dt}^2 \, \mathbb{E}\left[ \left( \frac{\partial \mathcal{L}_i}{\partial u} - \frac{\partial}{\partial t} \frac{\partial \mathcal{L}_i}{\partial u_t} \right)^2\right] = \mathbb{E}\left[\left(\frac{\partial \sigma_{i}}{\partial u} {d \mathcal{W}_t} \right)^2\right] .
\end{equation}
By utilizing the concept of quadratic variation from It\^{o} calculus, i.e., $\mathbb{E}[(d\mathcal{W})^2] = dt$, $\lim_{dt \to 0}(dt)^2 = 0$, $\mathbb{E}[d\mathcal{W}dt] = 0$ and assuming that the spatial and temporal white noise processes are uncorrelated to each other, Eq. \eqref{var_SL} can be rewritten as,
\begin{equation} \label{var_SL2}
    \lim_{dt\to 0} \, dt \, \mathbb{E}\left[ \left( \frac{\partial \mathcal{L}_i}{\partial u} - \frac{\partial}{\partial t} \frac{\partial \mathcal{L}_i}{\partial u_t} \right)^2\right] = \mathbb{E}\left[ \left(\frac{\partial \sigma_{i}}{\partial u}\right)^2 \right] ,
\end{equation}
where the squared white noise volatility in the right-hand side of Eq. \eqref{var_SL2} can assumed to be a linear combination of basis functions with coefficient $\boldsymbol{\beta}_i$, for $i=1, \ldots,n$. Similar to the Lagrangian density, we incorporate these basis functions into the library ${\mathcal{G}}\left(u, u_t, u_x\right)$, which can be different from the one used for the discovery of Lagrangian as the dependency of the diffusion is limited to a finite type of functions of states. 
To that end, we rewrite the Eq. \eqref{var_SL2} in terms of the library ${\mathcal{G}}\left(u, u_t, u_x\right)$ as,
\begin{equation} \label{var_SL3}
    \lim_{{dt\to 0}} \, dt \, \mathbb{E}\left[ \left( \frac{\partial \mathcal{L}_i}{\partial u} - \frac{\partial}{\partial t} \frac{\partial \mathcal{L}_i}{\partial u_t} \right)^2\right] = \mathbb{E}\left[ {\mathcal{G}}\left(u, u_t, u_x\right) \right]{\boldsymbol{\beta}}_i^{\top} .
\end{equation}
An example of the diffusion library employed in this work for uncovering the volatility functions of stochastic force is,
\begin{equation}
    \begin{aligned}
        \mathcal{G}\left(u, u_t, u_x\right) = {\left[\begin{array}{llllllllll}
        u & u_t & u_x & \mathrm{P}^{2}(u) & \mathrm{P}^{2}(u_t) & \mathrm{P}^{2}(u_x) & \operatorname{sin}(u) & \operatorname{cos}(u_t) & u\vert u\vert \\
        \end{array}\right], }
    \end{aligned}
\end{equation}
where $\mathrm{P}^{2}(\cdot)$ denotes the second-order polynomial function, sin($\cdot$) denotes the sine function, and cos($\cdot$) denotes the cosine function. 
Here, the problem boils down to identifying the coefficient vector $\boldsymbol{\beta}_i$, which can be efficiently obtained using the aforementioned sequential threshold least-square regression method \cite{tripura2024discovering}. 
Employing multiple measurements to account for the inherent stochastic variability in the excitation process, our aim is to obtain a least-squares solution for the coefficient ${\bm{\beta}}_{i}$. The estimation of ${\bm{\beta}}_i$ proceeds through analogous steps to the estimation of ${\bm{C}}_i$, as delineated in section \ref{sub1}.
However, upon sparse regression, a squared root operation would be required to be performed on the discovered expression from $\mathbb{E}\left[ {\mathcal{G}}\left(u, u_t, u_x\right) \right]{\boldsymbol{\beta}}_i^{\top}$ to correctly distill the terms in ${\partial \sigma_{i}}/{\partial u}$. 
Algorithm \ref{algo2} provides the implementation of the proposed framework for the discovery of volatility of stochastic excitation. 
\begin{algorithm}[t!]
    \small
    \caption{Discovery of diffusion term for stochastically excited system}
	\begin{algorithmic}
    \Require{Observed responses: $\{u_i(t), (u_x(t))_{i}, (u_t(t))_{i}\}^{N}_{j=1} \in \mathbb{R}^{N_t \times n \times N}$, the discovered Lagrangian density $\mathcal{L}_{i}$, and sparsification constant $\lambda \in \mathbb{R}$.} 
    \State {Select $m$ candidates for generalized stochastic force.}
    \State {Construct a symbolic library $\mathcal{G}\left( u, u_x, u_t \right) \in \mathbb{R}^{m}$.} \Comment{Can be different from Algorithm \ref{algo1}}
        \For{$i = 1$ : $n$}\Comment{Iteration for all particles}
        \State{Compute the derivatives of Lagrangian $\mathcal{L}_{i}$ over ensemble: $\{ \frac{\partial \mathcal{L}_i}{\partial u}, \frac{\partial}{\partial t} \bigl( \frac{\partial \mathcal{L}_i}{\partial u_t} \bigr) \}_{j=1}^{N}$. }
        \State {Compute basis functions of diffusion library over the ensemble: $\{ \mathcal{G}\left( u, u_x, u_t \right) \}_{j=1}^{N}$.}
        \State {Obtain target of regression: $dt \, \mathbb{E}\left[ \left( \frac{\partial \mathcal{L}_i}{\partial u} - \frac{\partial}{\partial t} \frac{\partial \mathcal{L}_i}{\partial u_t} \right)^2\right]$.}
        \State {Construct the regression equation.} \Comment{Refer Eq. (\ref{var_SL3})} 
        \State {Do sequentially threshold least-squares to obtain the coefficient ${\boldsymbol{\beta}}_i$.} \Comment{Refer Algorithm \ref{algo}} 
        \State{Do square root operation to uncover the volatility function: $\frac{\partial \sigma_{i}}{\partial u} = \left( \mathbb{E}\left[ {\mathcal{G}}\left(u, u_t, u_x\right) \right]{\boldsymbol{\beta}}_i^{\top} \right)^{1/2}$.}
        \EndFor
    \Ensure{Discovered volatility of stochastic excitation in terms of generalized force.}
	\end{algorithmic}
	\label{algo2}
\end{algorithm}

The discovery of the Lagrangian density and volatility function of diffusion provides an understanding of complex stochastic dynamical systems. The incorporation of sparse regression for Lagrangian discovery, rooted in a symbolic library encompassing various energy potentials, represents a significant methodological stride. In the next section, through numerical experiments, we will show that the proposed framework not only ensures the interpretability in Lagrangian modeling but also demonstrates versatility in addressing stochastically driven dynamical systems from Stochastic Differential Equations (SDEs) to Stochastic Partial Differential Equations (SPDEs).

\section{Numerical experiments}\label{sec:numerical} 
We consider six synthetic examples, four stochastic discrete systems and two stochastic continuous systems, to test the performance of the proposed methods. The examples are the stochastic harmonic oscillator, stochastic pendulum, stochastic Duffing oscillator, stochastic 3DOF structure, stochastic wave equation, and the stochastic Euler-Bernoulli Beam. The training data for the discrete systems are generated by solving the corresponding SDEs numerically using the Order 1.5 Strong Taylor Scheme \cite{kloeden1992higher}. The training data for the continuous systems are generated by solving the corresponding SPDEs numerically using the Order 0.5 Euler-Maruyama scheme \cite{kloeden1992higher}. 
For each example, we generate an ensemble of state trajectories, which includes the displacement and velocity responses. Each trajectory is simulated up to a time $t_f$ for different Wiener processes $\mathcal{W}_t$. The size of the ensemble and $t_f$ are mentioned in the corresponding sections. Further, we assume that only the ensemble of state trajectories is measurable and the input stochastic force vector, i.e., the Wiener process, is unknown to us.
Assuming that the length of each realization in the ensemble is $N_t$ (the length of time series data), we evaluate the candidate library and obtain a $\mathbb{R}^{N_t \times m}$ dimensional library $\mathbf{D} \in \mathbb{R}^{N_t \times m}$. 
For computing the time and spatial derivatives involved in the features of the design library $\mathbf{D}$, we use the second-order central finite difference,
\begin{equation} 
    \frac{d \omega}{d \eta} = \frac{\omega(\eta + \Delta \eta) -2\omega(\eta) + \omega(\eta - \Delta \eta)}{\Delta \eta^2} , 
\end{equation}
where $\omega$ is a function of $\eta$ and $\Delta \eta > 0$ is the increment in $\eta$. We also make use of the Lagrange's interpolation to calculate the derivatives, given as,
\begin{equation} 
    \begin{aligned}
        \frac{d \omega}{d \eta} &= \frac{\left(2 sh_1 - h_2 \right)}{h_1 (h_1 + h_2)}\omega_i + \frac{\left(2 s +1 \right)h_1 - h_2}{h_1 h_2}\omega_{i+1} + \frac{ \left(2s + 1 \right)h_1 }{h_2 (h_1 + h_2)}\omega_{i+2}, \\
        \frac{d^2 \omega}{d \eta^2} &= \frac{2 \left(h_2 \omega_i - (h_1+h_2)\omega_{i+1} + h_1 \omega_{i+2} \right)}{h_1 h_2 (h_1 + h_2)} , 
    \end{aligned}
\end{equation}
where $h_1=\omega_{i+1}-\omega_{i}$ and $h_2=\omega_{i+2}-\omega_{i+1}$. We can find the derivative values at three neighborhood grid points by substituting $s$ equal to -1, 0, and $h_2/h_1$. The above formula is also useful for computing derivative in irregular observations, and when $h_1=h_2=h$, the above formula becomes a three-point forward finite difference. 
We evaluate the library $\mathbf{D}$ on each realization in the ensemble, which yields the set $\{\mathbf{D} \in \mathbb{R}^{N_t \times m} \}_{j=1}^{N}$ with $N$ denoting the ensemble size. Thereafter, the mean of the ensemble of libraries is taken for obtaining the mean library $\mathbb{E}[\mathbf{D}] \in \mathbb{R}^{N_t \times m}$. The final sparse regression is performed on the mean design library.

\subsection{Stochastic Harmonic oscillator}\label{example_harmonic}
\begin{figure}[ht!]
    \centering
    \includegraphics[width=\textwidth]{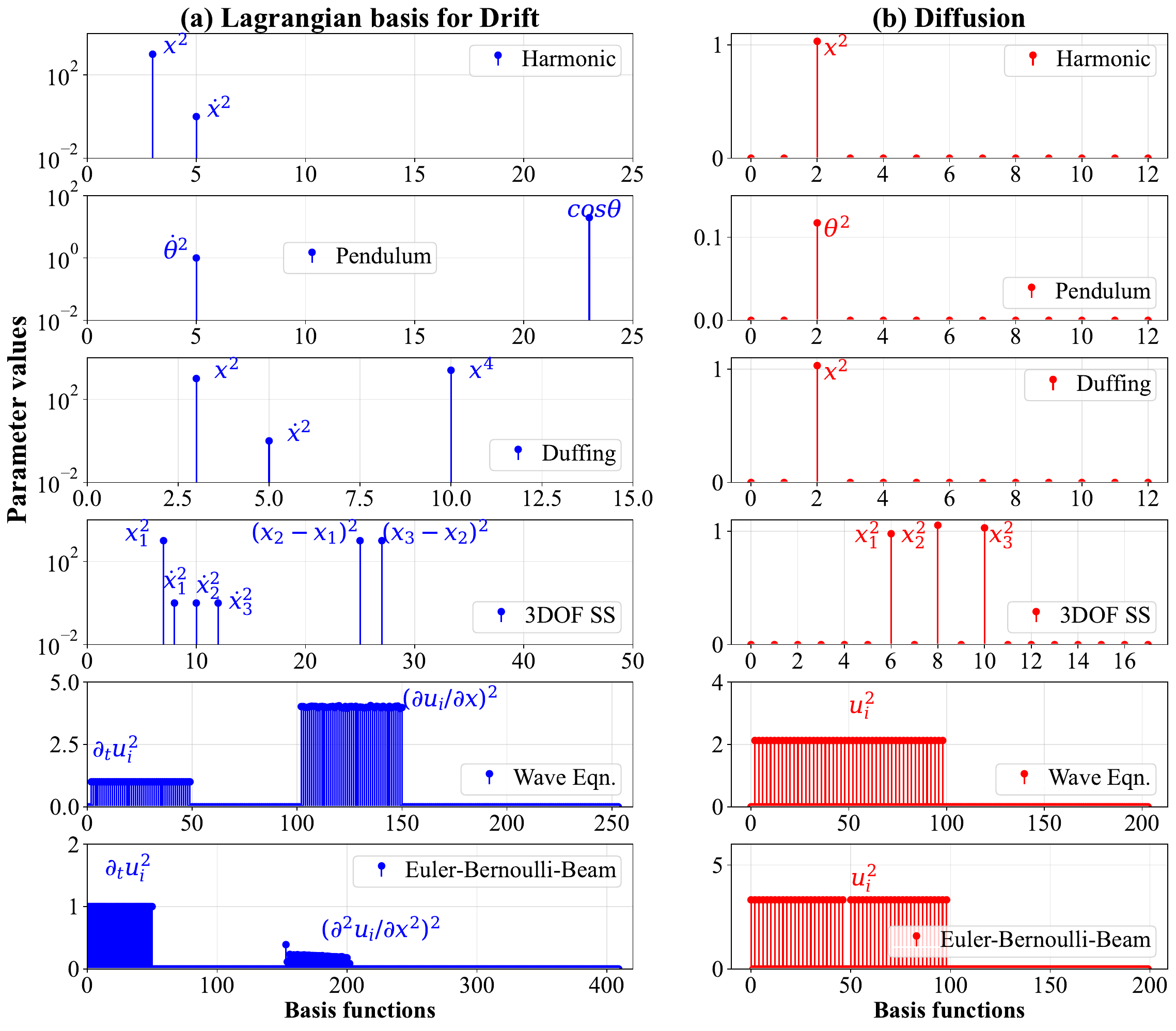}
    \caption{Discovery of the Lagrangian basis functions for the example SDEs and SPDEs. The identified basis functions are shown in the x-axis, and the associated parameters are shown in the y-axis. (a) The plots on the left side show the identified basis functions of deterministic Lagrangian density. The library matrix $\mathbf{D}^{\text{f}}$ for the deterministic Lagrangian density contains 25 basis functions for the stochastic Harmonic oscillator and stochastic pendulum, 15 basis functions for the stochastic Duffing oscillator, 50 for the stochastic 3DOF structural system, 254 for the stochastic wave equation, and 421 for the stochastic Euler-Bernoulli Beam. (b) The plots on the right show the potential energy basis functions associated with the diffusion process. In this case, the library matrix $\mathbf{D}^{\text{g}}$ contains 12 basis functions for the stochastic harmonic, pendulum, and Duffing oscillator, 17 for the 3DOF system, 204 for the stochastic wave equation, and 200 for the stochastic Euler-Bernoulli Beam equation.}
    \label{fig:basis}
\end{figure}

In the first example, we consider the stochastic harmonic oscillator, whose governing equation of motion is given as,
\begin{equation}\label{eq:harmonic}
    m\ddot{X}(t) + kX(t) = \sigma \dot{\mathcal{W}}(t); \;\; X(t=0)=X_0, \; t \in [0,t_f],
\end{equation}
where $m \in \mathbb{R}$ and $k \in \mathbb{R}$ are the mass and spring stiffness of the oscillator. $\dot{\mathcal{W}} \sim \mathcal{N}(0,1)$ is the white noise, expressed as the derivative of the Wiener process $\mathcal{W} \sim \mathcal{N}(0,t)$ and $\beta \in \mathbb{R}$ is the strength of the white noise. To solve the system, we introduce the statespace variables $\{Y_1=X, Y_2=\dot{X}\}$ and write the SDE representation of the above equation as,
\begin{equation}\label{eq:sde_harmonic}
    \begin{aligned}
        dY_1(t) &= Y_2(t)dt, \\
        dY_2(t) &= -\frac{k}{m} Y_1(t)dt + \frac{\sigma}{m} d\mathcal{W}(t) .
    \end{aligned}
\end{equation}
The system properties are taken as $m$=1kg, $k$=1000N/m, and $\sigma$=1. For generating training data, we consider an initial condition of $\{Y_1, Y_2\}(0)$=\{0.5,0\} and solve the system for $t_f=1$s using a sampling frequency of 10000Hz. We generate $N_e$=200 realizations of the state trajectories (i.e., $\{Y_1(t), Y_2(t)\}_{i=1}^200$) for different $\mathcal{W}(t)$. The true Lagrangian density function of the oscillator is provided in Table \ref{table:lagrange}. 
The identified basis functions from the sparse regression are illustrated in Fig. \ref{fig:basis}(a) and \ref{fig:basis}(b). As depicted in Fig. \ref{fig:basis}(a), the proposed framework identifies $X^2$ and $\dot{X}^2$ as the basis functions of Lagrangian density and $X^2$ as the basis function for Wiener potential. With the identified parameters, the identified Lagrangian density function is given as $\mathcal{L}^{*}_{s}$=$0.5\dot{X}^2 - 500.06X^2$, whereas the true Lagrangian is $\mathcal{L}_{s}$=$0.5\dot{X}^2 - 500X^2$. The Wiener potential is identified as $1.03 X$, whereas the true value is $X$. Here, note that the Wiener potential is obtained by taking the square root of the identified squared Wiener potential function $X^2$. These expressions suggest that the proposed framework is able to identify the Lagrangian density function and the Winer potential of the diffusion functions almost exactly.

\begin{table}[ht!]
	\centering
	\caption{Lagrangian description of the undertaken examples.}
	\label{table:lagrange}
    \small
	\begin{tabular}{lll}
		\toprule
		{System} & {Lagrangian density function} & {Wiener potential} \\
		\midrule
		Stoch. Harmonic oscillator & $\mathcal{L}_{s}$=$0.5\dot{X}^2 - 0.5{k}{m^{-1}}X^2$ & $\sigma x / m$ \\
		Stoch. pendulum & $\mathcal{L}_{s}$=$0.5\dot{\theta}^2 - {g}{l^{-1}}{\rm{cos}}\theta$ & $\sigma \theta / m$ \\
		Stoch. Duffing oscillator & $\mathcal{L}_{s}$=$0.5\dot{X}^2 - 0.5{k}X^2 - 0.25{\alpha}X^4$ & $\sigma x$ \\
		Stoch. 3DOF system & $\mathcal{L}_{s}$=$0.5\sum_{k=1}^{3} m_k \dot{X}_k^2 - 0.5\left( k_1 X_1^2 + k_2(X_2-X_1)^2 + k_3(X_3-X_2)^2 \right)$ & $\sum_{i=1}^3 \sigma_i x_i / m_i$ \\
		Stoch. wave equation & $\mathcal{L}_{s}$=$0.5\sum_{i}{(\partial_t u_i(x,t))^2} - c^2 \sum_{i}{(\partial_x u_i(x,t))^2}$ & $\sum_{i} \sigma_i c^2 u_i$ \\
		Stoch. Euler-Bernoulli Beam & $\mathcal{L}_{s}$=$0.5\sum_{i}{(\partial_t u_i(x,t))^2} - 0.5{EI}{(\rho A)^{-1}} {(\partial_{xx} u_i(x,t))^2}$ & $\sum_{i=1}^3 \sigma_i \rho A u_i$ \\
		\bottomrule
	\end{tabular}
\end{table}

\subsection{Stochastic pendulum}\label{example_3}
Next, we consider the vibration of the stochastic pendulum under the action of gravity and external stochastic force. The governing equation of motion of a simple pendulum is given as follows,
\begin{equation}
    m l^2 \ddot{\theta}(t) + mgl \; {\rm{sin}} \theta(t) = \sigma \dot{\mathcal{W}}(t); \;\; \theta(t=0)=\theta_0, \\; t \in [0,t_f] ,
\end{equation}
where $m \in \mathbb{R}$ is the mass of the pendulum, $l \in \mathbb{R}$ is the length of the suspension, and $g \in \mathbb{R}$ is the gravitational acceleration constant. Like the previous case, $\mathcal{W}$ represents the stochastic Wiener noise with strength $\sigma \in \mathbb{R}$. We use $m$=1kg, $l$=1m, and $\sigma$=0.1 to simulate training data. We write the corresponding SDE, similar to Eq. \eqref{eq:sde_harmonic}, and use the Order 1.5 Strong Taylor scheme to solve the system. We use a sampling frequency of 2000Hz to simulate the system for a total duration of $t_f$=5s and generate $N_e$=200 realizations. The initial angular displacement and velocity are considered as $\{\theta, \dot{\theta}\}(0)$=\{0.9,0\}. 
Table \ref{table:lagrange} provides the underlying Lagrangian description. 
The identification results for the Lagrangian density and the Wiener potential are presented in Fig. \ref{fig:basis}(a) and \ref{fig:basis}(b). Out of the 25 candidate functions, only $\dot{\theta}^2$ and $\operatorname{cos}\theta$ are identified as the basis functions of Lagrangian. With the identified parameters, the identified Lagrangian density is given as $\mathcal{L}^{*}_{s}$=$0.5\dot{\theta}^2 - 9.80{\operatorname{cos}}\theta$. By substituting the appropriate values, we can confirm from the Table \ref{table:lagrange} that the actual Lagrangian is $\mathcal{L}_{s}$=$0.5\dot{\theta}^2 - 9.81{\operatorname{cos}}\theta$. Similarly, after taking the square root operation, the identified and true Wiener potentials are $0.11 \theta$ and $0.10 \theta$. We observe a close match in both the Lagrangian and Wiener potential of diffusion terms, indicating the proposed framework's robustness.

\subsection{Stochastic Duffing oscillator}
Next, we consider the nonlinear stochastic Duffing oscillator, whose governing equation of motion is given as follows,
\begin{equation}\label{eq:duffing}
    \ddot{X}(t) + kX(t) + \alpha X^3(t) = \sigma \dot{\mathcal{W}}(t); \;\; X(t=0)=X_0, \; t \in [0,t_f],
\end{equation}
where $k \in \mathbb{R}$ is the stiffness of the linear spring, and $\alpha \in \mathbb{R}$ is the nonlinear dissipative spring constant. ${\mathcal{W}}(t)$ here denotes the usual stochastic Wiener noise with strength $\sigma \in \mathbb{R}$. An SDE representation of the above equation, similar to the Eq. \eqref{eq:sde_harmonic}, is derived for solving the system. 
We consider the spring constants as $k$=1000N/m and $\alpha$=5000. The noise strength and the initial conditions are taken as $\sigma$=1 and $\{X, \dot{X}\}(0)$=\{0.4,0\}. A total of $N_e$=200 realizations are generated by solving the system for $t_f$=1s using $\Delta t$=0.0001s.  
Table \ref{table:lagrange} provides the underlying Lagrangian description of the stochastic Duffing oscillator. By substituting the appropriate parameter values, we obtain true Lagrangian density as $\mathcal{L}_{s}$=$0.5\dot{X}^2 - 500X^2 - 1250X^4$ and the Wiener potential as $X$.
The corresponding identification results for the Lagrangian and Wiener potential are shown in Fig. \ref{fig:basis}. A total of 15 candidate basis functions are considered for discovering the Lagrangian density, and 12 basis functions are considered for discovering Winener potential. Out of the 15 basis functions, we identify $X^2$, $\dot{X}^2$, and $X^4$ for Lagrangian density, and out of 12 basis functions, we obtain $X^2$ for Winer potential. The identified Lagrangian density function is $\mathcal{L}^{*}_{s}$=$0.5\dot{X}^2 - 499.67X^2 - 1251.15X^4$ and the Wiener potential is $1.03 X$. The results suggest that the identified interpretable expressions are very close to the ground truth expressions.

\subsection{Stochastic 3DOF oscillator}\label{example_4}
Here, we consider a three degrees of freedom (DOF) stochastic undamped structural system. We denote the mass and stiffness of each DOF as $m_i$ and $k_i$. The equation of motion for the three DOF system follows,
\begin{equation}\label{eq:mdof}
    \mathbf{M} \ddot{\bm{X}}(t) + \mathbf{K} \bm{X}(t) = \mathbf{S} \dot{\bm{\mathcal{W}}}(t), \;\; \bm{X}(0) = \bm{X}_0, \; t\in [0,t_f],
\end{equation}
where the mass matrix $\mathbf{M} \in \mathbb{R}^{3 \times 3}$, stiffness matrix $\mathbf{K} \in \mathbb{R}^{3 \times 3}$, the volatility constant matrix $\mathbf{S} \in \mathbb{R}^{3 \times 3}$, and the Wiener noise vector $\bm{\mathcal{W}} \in \mathbb{R}^3$ are given as,
\begin{equation}
    \mathbf{M} = \left[ { \begin{array}{ccc}
         m_1 & 0 & 0 \\
         0 & m_2 & 0 \\
         0 & 0 & m_3
    \end{array} } \right], \;
    \mathbf{K} = \left[ { \begin{array}{ccc}
         k_1+k_2 & -k_2 & 0 \\
         -k_2 & k_2+k_3 & -k_3 \\
         0 & -k_3 & k_3
    \end{array} } \right], \;
    \mathbf{S} = \left[ { \begin{array}{ccc}
         \sigma & 0 & 0 \\
         0 & \sigma & 0 \\
         0 & 0 & \sigma
    \end{array} } \right], \;
    \bm{\mathcal{W}} = \left[ { \begin{array}{c}
         \dot{\mathcal{W}}_1, \\
         \dot{\mathcal{W}}_2, \\
         \dot{\mathcal{W}}_3
    \end{array} }
    \right]
\end{equation}
To solve the system using Order 1.5 Strong Taylor scheme, we write the first-order It\^{o} SDEs by introducing the statespace $\{X_1,\dot{X}_1,X_2,\dot{X}_2,X_3,\dot{X}_3 \}=\{Y_1,Y_2,Y_3,Y_4,Y_5,Y_6\}$ as follows,
\begin{equation}\label{eq:sde_3dof}
    \begin{aligned}
        d \bm{Y}(t) = \bm{f}(\bm{Y},t)dt + \bm{\sigma} d\bm{\mathcal{W}}(t),
    \end{aligned}
\end{equation}
where $\bm{f}(\bm{Y},t) \in \mathbb{R}^{6}$ and $\bm{\sigma} \in \mathbb{R}^{6 \times 3}$ are drift and volatility matrix of the stochastic 3DOF system. We use $m_i$=10kg, $k_i$=10000N/m, and $\sigma_i$=1 for generating training data. The initial conditions are taken as $\bm{X} = \{0.25,0,0.5,0,0,0\}$. With a sampling frequency of 10000Hz, we generate $N_e$=200 realizations of system states; each simulated up to $t_f$=1s. The actual Lagrangian density function is described in Table \ref{table:lagrange}. Upon substitution of actual parameter values, we obtain $\mathcal{L}_{s}$=$0.5\sum_{k=1}^{3} \dot{X}_k^2 - (1000X_1^2 + 1000(X_2-X_1)^2 + 1000(X_3-X_2)^2 )$ and the Wiener potential as $\sum_{i=1}^3 0.01 X_i$.
The identified basis functions for the Lagrangian density and Wiener potential are given in Fig. \ref{fig:basis}. We observe the identified basis functions for the Lagrangian density as $X_1^2$, $(X_2-X_1)^2$, and $(X_3-X_2)^2$, and for the Wiener potential as $X_1^2$, $X_2^2$, and $X_3^2$. The final identified Lagrangian expression is $\mathcal{L}^{*}_{s}$=$0.5\sum_{k=1}^{3} \dot{X}_k^2 - (999.99X_1^2 + 1000.04(X_2-X_1)^2 + 999.96(X_3-X_2)^2 )$ and the Wiener potential is $0.97X_1 + 1.04X_2 + 1.01X_3$. The expressions show that the identified interoperable expressions are sufficiently accurate as the ground truth.

\subsection{Stochastic Wave equation}\label{example_6}
As a first continuous system, we consider the stochastic wave equation, often used to represent the traveling waves in the vibration of fluids and solids. The stochastic wave equation follows,
\begin{equation}
    \begin{aligned}
        {\partial_{tt} u(X,t)} &= c^2 {\partial_{xx} u(X,t)} + \sigma \partial_{t} \mathcal{W}(t), \;\; x\in [0,L], \; t\in [0,t_f], \\
        u(x,t) &= u(L,t) = 0, \\
        u(x,0) &= \operatorname{cos}(2\pi L),
    \end{aligned}
\end{equation}
where $c=\sqrt{s / \rho}$ is the velocity of travelling wave, $s \in \mathbb{R}$ is the shearing modulus, and $\rho \in \mathbb{R}$ is the density of the material. $\mathcal{W}(x,t)$ is the space-time Wiener noise whose time derivative provides the space-time white noise $\partial_{t} \mathcal{W}$, and $\sigma$ is the strength of the noise. For generating training data, we use $c$=2m/s and $\sigma$=2.
A total of $N_e$=30 realizations of systems response are generated by solving the above equation up to $t_f$=1s using $\Delta t$=0.0001s and $\Delta x$=0.01s. 
The associated Lagrangian is provided in Table \ref{table:lagrange}, whose final expression after the substitution of parameter values is $\mathcal{L}_{s}$=$0.5\sum_{i}{(\partial_t u_i(x,t))^2} - 4 \sum_{i}{(\partial_x u_i(x,t))^2}$ and the Wiener potential is $\sum_{i} 8 u_i$.
The identified basis functions for the Lagrangian and Wiener potential at the discretized spatial locations are illustrated in Fig. \ref{fig:basis}. In this example, for discovering the Lagrangian density of wave, we use 254 basis functions in the library, out of which only $(\partial_t u_i)^2$ and $(\partial_x u_i)^2$ are discovered at the spatial locations. Similarly, 204 basis functions are used in the library of Wiener potential, out of which only $u_i^2$ at each location is identified. The final identified Lagrangian density is given as $\mathcal{L}^{*}_{s}$=$0.5\sum_{i}{(\partial_t u_i(x,t))^2} - 4.0013 \sum_{i}{(\partial_x u_i(x,t))^2}$ and the Wiener potential is given as $\sum_{i} 8.5227 u_i$. 
Compared to the previous examples, although the identified Lagrangian function has a very small identification error, we observe a relatively higher discrepancy in the identified Wiener potential, amounting to 6.53\%.

\subsection{Stochastic Euler-Bernoulli Beam}\label{example_7}
As a second continuous system, we consider the stochastic Euler-Bernoulli beam with a homogeneous and uniform cross-section. The following equation of motion expresses the stochastic motion of the Euler-Bernoulli beam,
\begin{equation}
    \begin{aligned}
        {\partial_{tt} u(x,t)} &= c {\partial_{xxxx} u(x,t)} + \sigma \partial_t {\mathcal{W}}(x,t), \;\; x \in [0,L], \; t\in [0, t_f], \\
        u(0,t) &= \partial_x u(0,t) = \partial_{xx}u(L,t) = \partial_{xxx}u(L,t) = 0 , \\
        u(x,0) &= (\operatorname{cosh} \psi x - \operatorname{cos} \psi x) + \frac{(\operatorname{cos} \psi L + \operatorname{cosh} \psi L)}{(\operatorname{sin} \psi L + \operatorname{sinh} \psi L)}(\operatorname{sin} \psi x - \operatorname{sinh} \psi x) ,
    \end{aligned} 
\end{equation}
where $c= {EI}/{\rho A}$, $E \in \mathbb{R}$ is the modulus of elasticity of the beam, $I \in \mathbb{R}$ is the moment-of-inertia, $\rho \in \mathbb{R}$ is the density of the beam, and $A \in \mathbb{R}$ is the cross-section of the beam. Here $\mathcal{W}(x,t)$ is the space-time Wiener noise, and $\sigma$ is the strength of the noise. For synthetic data generation, we use the following system properties: $E$=$2\times10^10$N/m$^2$, $A$=$0.02\times 0.001$m$^2$, $\rho$=8050kg/m$^3$, and $L$=1. The beam's natural frequency is taken as $\psi$=0.596864$\pi$, and the noise strength $\sigma$ is considered as 20.
The beam is solved for $t_f$=2s using $\Delta t$=0.0001s and $\Delta x$=0.01. A total of $N_e$=20 realizations are simulated to obtain the mean statistics of the beam responses. The underlying Lagrangian is provided in Table \ref{table:lagrange}, whose exact description for the considered parameters is $\mathcal{L}_{s}$=$0.5\sum_{i}{(\partial_t u(x,t))^2} - 0.1035 {(\partial^2_x u(x,t))^2}$. For the Wiener potential, the exact description is $\sum_{i=1}^3 3.22 u_i$.
The identified basis functions are illustrated in Fig. \ref{fig:basis}. Here, the library for Lagrangian discovery contains 421 candidate basis functions, and the library for Wiener basis consists of 200 basis functions. In the Lagrangian discovery, we identify the basis functions $(\partial_t u_i)^2$ and $(\partial_{xx} u_i)^2$ at the spatial locations. In the Wiener potential library, we identify $u_i^2$ at the spatial locations. The final identified Lagrangian description is $\mathcal{L}^{*}_{s}$=$0.5\sum_{i}{(\partial_t u(x,t))^2} - 0.1024 {(\partial^2_x u(x,t))^2}$, and the identified Wiener potential description is $\sum_i 3.2216 u_i$. In both descriptions, we identify a relative error of 1.0628\% and 0.05\%, which again indicates the robustness of the proposed framework.

\subsection{Identification of SDEs from the discovered Lagrangian}
In this section, we discover the interpretable governing equations of motion of underlying dynamical systems by evaluating the Euler-Lagrangian equation. The discovered Lagrangian equations of motion are compared with the ground truth in Table \ref{table:eqm}. For quantitative evaluation, we estimate the percentage relative L$^2$ norm of the error between the true and identified system parameters as $100 \times \frac{\| \bm{\theta}-\bm{\theta}^*\|}{\|\bm{\theta}\|}$, where $\bm{\theta}$ and $\bm{\theta}^*$ represents the true and identified system parameters. We observe that the identification error is less than one percent in all the cases except the stochastic wave equation, where we observe a relative error of approximately 2.9\%. The ability to almost exactly approximate the true equations of motion using only one second of data is a novelty of the proposed framework. 
\begin{table}[b!]
	\centering
	\caption{Summary of the discovered equations of motion of example problems. The error represents the relative L$^2$ error between the parameters of the true and identified stochastic systems, including the white noise strength.}
	\label{table:eqm}
    \small
	\begin{tabular}{llll}
		\toprule
		{System} & \multicolumn{2}{l}{Equation of motion of underlying SDEs/PDEs} & Error (\%) \\
		\midrule
		Stoch. harmonic  & Truth: & $\ddot{X}(t) + 1000X(t) = \dot{W}(t)$ & -- \\
        & \textit{Discovered}: & $\ddot{X}(t) + 1000.12X(t) = 1.03\dot{W}(t)$ & 0.0116 \\
		Stoch. pendulum & Truth: & $\ddot{\theta}(t) + 9.81{\rm{sin}} \theta(t) = 0.10\dot{W}(t)$ & -- \\
        & \textit{Discovered}: & $\ddot{\theta}(t) + 9.80{\rm{sin}} \theta(t) = 0.11\dot{W}(t)$ & 0.1441 \\
		Stoch. Duffing  & Truth: & $\ddot{X}(t) + 1000X(t) + 2500 X^3(t) = \dot{W}(t)$ & -- \\
        & \textit{Discovered}: & $\ddot{X}(t) + 999.35X(t) + 2502.31 X^3(t) = 1.03\dot{W}(t)$ & 0.0891 \\
		Stoch. 3DOF  & Truth: & $\ddot{X}_1(t) + 1000X_1(t) - 1000(X_2-X_1)(t) = \dot{W}_1(t)$ & -- \\
        && $\ddot{X}_2(t) + 1000(X_2-X_1)(t) - 1000(X_3-X_2)(t) = \dot{W}_2(t)$ & -- \\
        && $\ddot{X}_3(t) + 1000(X_3-X_2)(t) = \dot{W}_3(t)$ & -- \\
        & \textit{Discovered}: & $\ddot{X}_1(t) + 999.99X_1(t) - 1000.04(X_2-X_1)(t) = 0.97\dot{W}_1(t)$ & 0.0036 \\
        && $\ddot{X}_2(t) + 1000.04(X_2-X_1)(t) - 999.96(X_3-X_2)(t) = 1.04\dot{W}_2(t)$ & 0.0048 \\
        && $\ddot{X}_3(t) + 999.96(X_3-X_2)(t) = 1.01\dot{W}_3(t)$ & 0.0041 \\
		Stoch. wave  & Truth: & ${\partial^2_t u(x,t)} = 4 {\partial^2_x u(x,t)} + 2\partial_t W(x,t)$ & -- \\
        & \textit{Discovered}: & ${\partial^2_t u(x,t)} = 4.0013 {\partial^2_x u(x,t)} + 2.13\partial_t W(x,t)$ & 2.9069 \\
		Stoch. beam  & Truth: & ${\partial^2_t u(x,t)} = 0.1035 {\partial^4_x u(x,t)} + 20\partial_t W(x,t)$ & -- \\
        & \textit{Discovered}: & ${\partial^2_t u(x,t)} = 0.1024 {\partial^4_x u(x,t)} + 20.01\partial_t W(x,t)$ & 0.0503 \\
		\bottomrule
	\end{tabular}
\end{table}
\begin{figure}[t!]
    \centering
    \includegraphics[width=\textwidth]{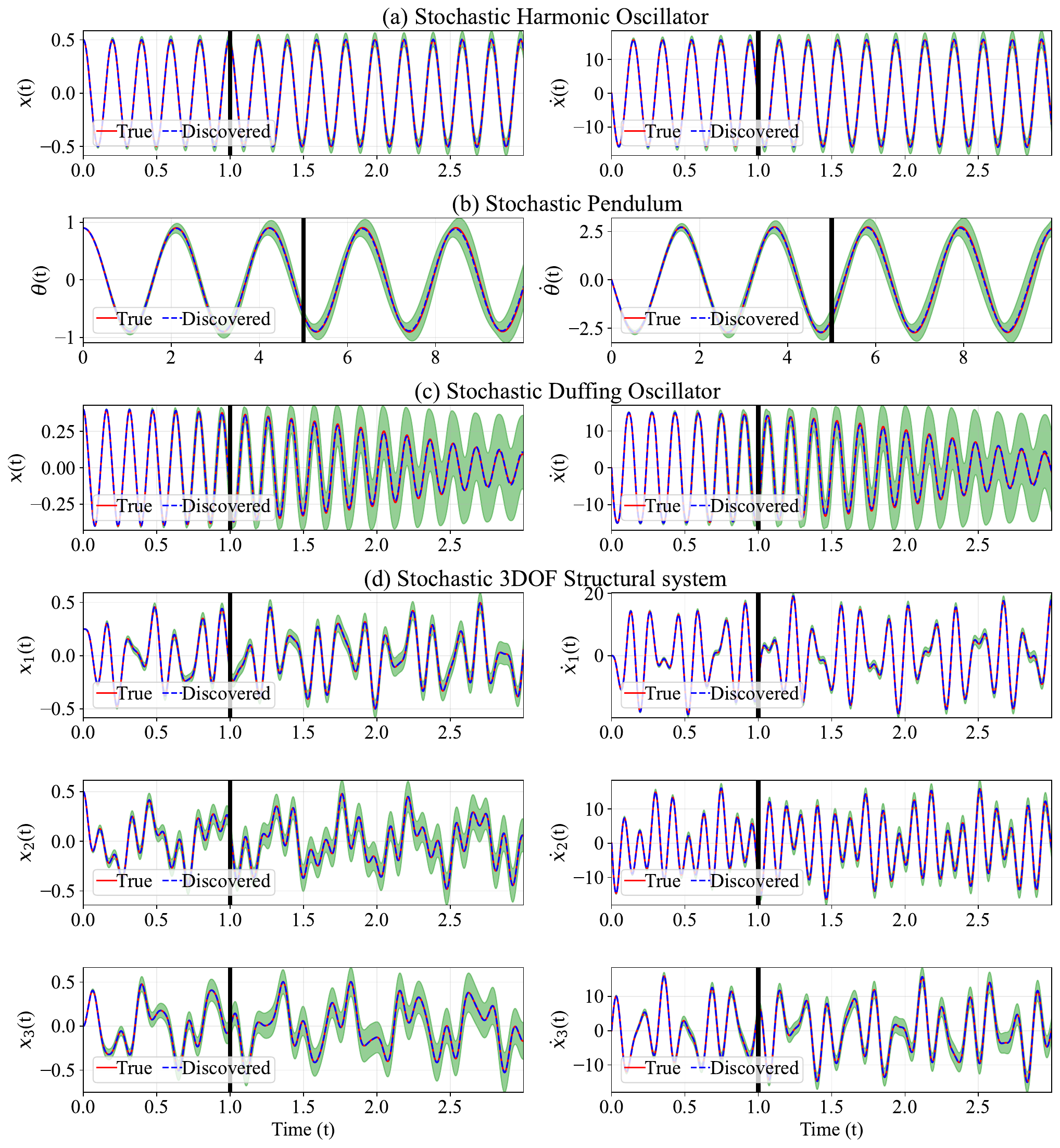}
    \caption{Time evolution of the responses of the true and identified stochastic differential equations. The solid line indicates the expected response of the true (red solid line) and discovered (dashed blue line) SDEs. The green shaded region indicates the uncertainty (two times the standard deviation) in the system response due to the stochastic Wiener process. The black vertical solid line separates the training (the left side) and prediction regions (the right side).}
    \label{fig:response_sde}
\end{figure}

We solve the identified SDEs using the strong Taylor order 1.5 and Euler-Maruyama schemes to visually compare the responses of true and identified systems. We simulate 200 realizations to estimate the mean and standard deviation of the stochastic responses. In Fig. \ref{fig:response_sde}, the mean responses of the identified discrete systems are compared with the mean responses of true systems. Also, the uncertainty in the identified systems, estimated using two times standard deviation, due to the stochastic Wiener process is illustrated in Fig. \ref{fig:response_sde}. The left side of the vertical black indicates the training regions and the right side indicates the prediction regions. In all the examples, we observe that the mean responses of the identified equations not only accurately match the responses of the true systems in the training regime but also almost exactly emulate the true responses in the prediction region. This indicates the generalizability of the identified equations beyond the training period. The proposed framework's ability to identify the equations of motion along with Lagrangian in a single framework is a second novel contribution to the literature.

\begin{figure}[t!]
    \centering
    \includegraphics[width=\textwidth]{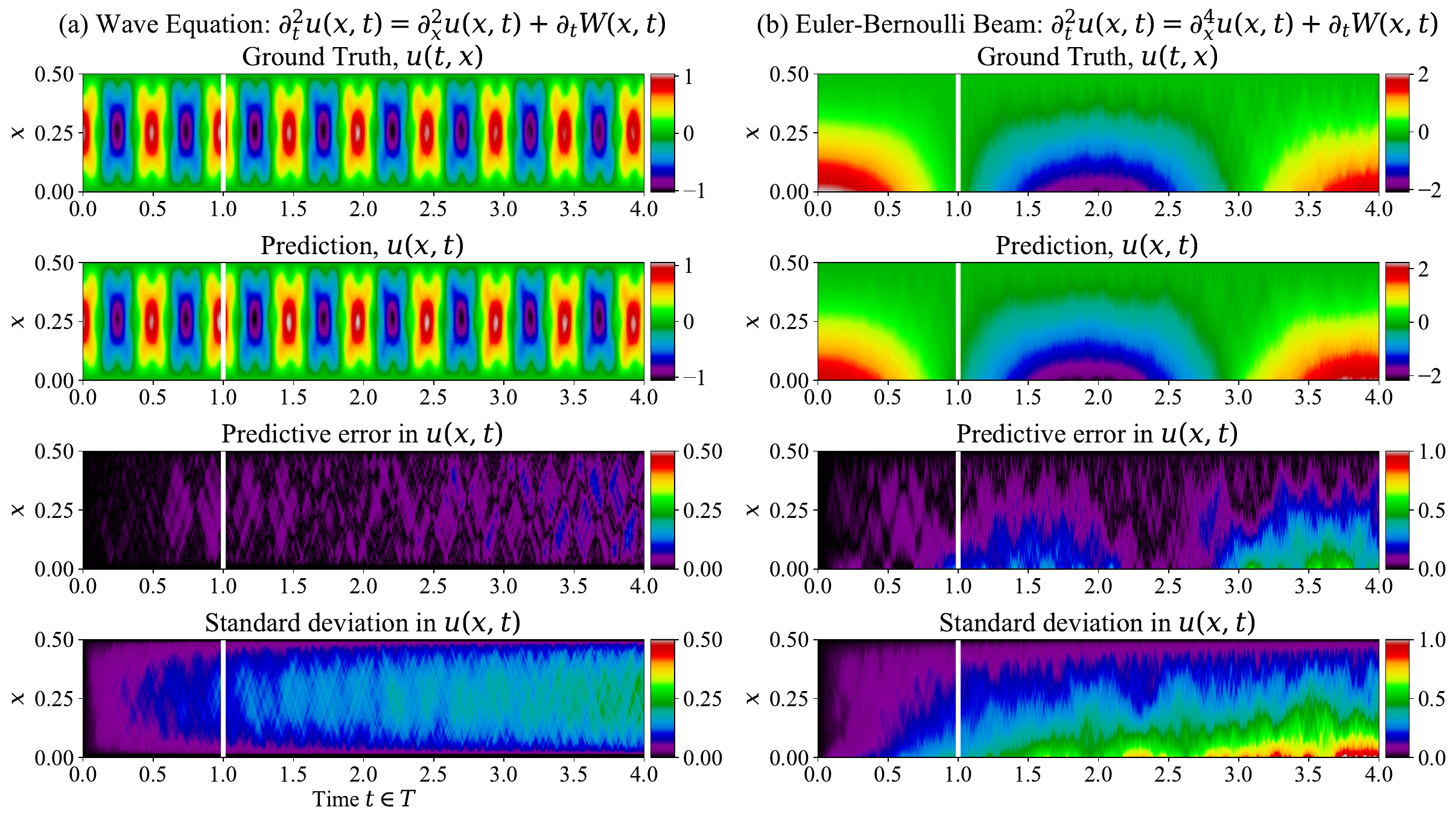}
    \caption{Time evolution of the response of the discovered wave equation and Euler-Bernoulli beam. The ground truth indicates the response of the true stochastic wave equation and stochastic Euler-Bernoulli Beam. The prediction indicates the response of the discovered stochastic wave equation and stochastic Euler-Bernoulli Beam. The predictive error indicates the absolute error between the ground truth and the predicted response. The standard deviation indicates the randomness in the system response due to the stochastic Wiener noise. The left side of the white vertical indicates the training duration, and the right side indicates the prediction duration.}
    \label{fig:response_pde}
\end{figure}

Similarly, in Fig. \ref{fig:response_pde}, the mean responses of the true and identified continuous systems are displayed. The left side of the vertical white line indicates the training period, and the right side indicates the prediction regime. Similar to the results of discrete systems, we observe that the identified equations almost exactly reconstruct the responses of the ground truth in the training regime. 
However, we observe a minor discrepancy in the responses of the identified systems with an increase in the prediction time. Consequently, the predictive error and standard deviation increase as the prediction time duration increases. 
In this work, the proposed algorithm utilizes only one second of training data, except for the stochastic pendulum example, to identify the governing equations of motion. Even with a small amount of data, the identified equations almost exactly emulate the ground truth responses. This can be attributed to the fact that the identified equations from the proposed framework follow the underlying conservation laws, thereby guaranteeing the true optimal dynamics. 
Nevertheless, the results reported here for the continuous systems can be further improved by using more training data. 

\subsection{Identifying Hamiltonian representations from the discovered Lagrangian}
\begin{figure}[ht!]
    \centering
    \includegraphics[width=0.8\textwidth]{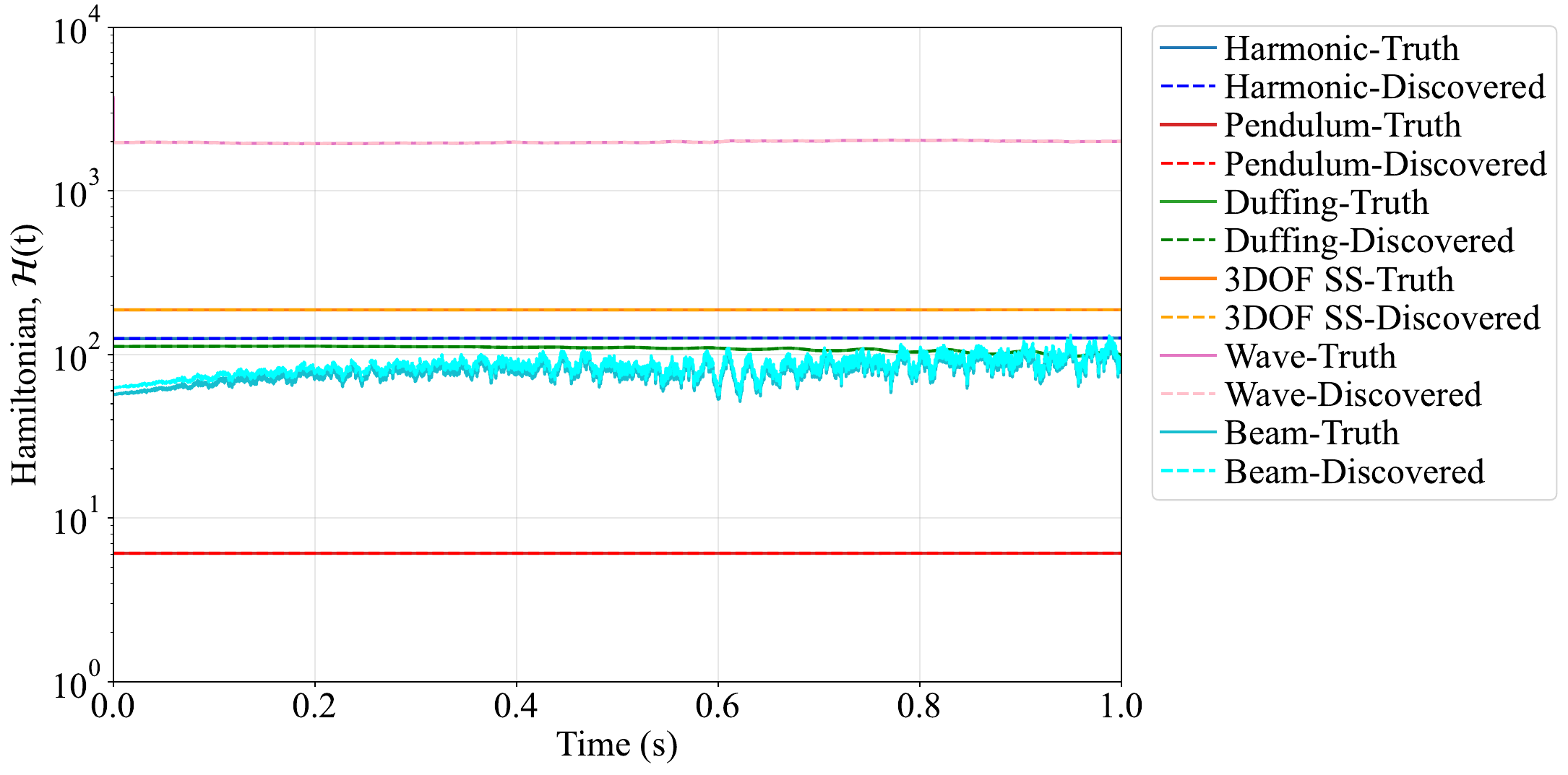}
    \caption{Time evolution of Hamiltonian of the example problems. The solid lines indicate the Hamiltonian trajectory of the true systems, and the dashlines indicate the Hamiltonian trajectory of the discovered systems.}
    \label{fig:hamiltonian}
\end{figure}

Hamiltonian is an important property of engineering systems, which provides fundamental knowledge about energy conservation in physical systems. Due to the importance of Hamiltonian, for e.g., control \cite{duong2021hamiltonian}, a significant number of Hamiltonian discovery frameworks are proposed in the literature. However, these frameworks are limited to deterministic systems only, i.e., where the input force information is known. In this context, the proposed framework not only eliminates the need for a separate Hamiltonian discovery framework but can also be applied to discover Hamiltonian for systems with stochastic excitation. Once the Lagrangian is discovered, the Hamiltonian can be discovered by applying the Legendre transformation to the discovered Lagrangian,
\begin{equation}
    \begin{aligned}
        \mathcal{H}_{s} &= \sum^n_{i=1} \frac{\partial \mathcal{L}_{s}}{\partial \dot{u}_{i}} \dot{u}_{i} - \mathcal{L}_{s} \\
        & = \sum^n_{i=1} \frac{\partial}{\partial \dot{u}_{i}} {\left( \sum_{j=1}^{n} \mathbb{E}[\mathbf{D}]\bm{C}_j \right)} \dot{u}_{i} - {\left( \sum_{j=1}^{n} \mathbb{E}[\mathbf{D}]\bm{C}_j \right)} .
    \end{aligned}
\end{equation}
The interpretable expressions of the identified Hamiltonians are provided in Table \ref{table:hamiltonian}. The identified parameters are presented in the brackets. By evaluating the true expressions of Hamiltonian from Table \ref{table:lagrange} and comparing the ground truth with identified expressions, one can observe that the identifying Hamiltonian expressions almost exactly match the ground truth. 
The time evolution of the Hamiltonian for the undertaken examples is displayed in Fig. \ref{fig:hamiltonian}. The solid line represents the Hamiltonian of the true system, and the dashed line indicates the Hamiltonian values from identified expressions in Table \ref{table:hamiltonian}. It is evident that the Hamiltonian values from the identified equation almost completely overlap the ground truth values, indicating the accuracy of the identified systems. The ability to accurately identify the Hamiltonian along with Lagrangian and equations of motion in a single framework from stochastic observation is another novelty of the proposed framework. 
\begin{table}[t!]
	\centering
	\caption{Summary of the discovered Hamiltonian of example problems}
	\label{table:hamiltonian}
    \small
	\begin{tabular}{ll}
		\toprule
		{System} & {Hamiltonian description} \\
		\midrule
		Stochastic Harmonic oscillator & $\mathcal{H}_{s}$=$\frac{1}{2}\dot{X}^2 + \frac{1}{2}\{\frac{k}{m}=1000.12\}X^2$ \\
		Stochastic pendulum & $\mathcal{H}_{s}$=$\frac{1}{2}\dot{\theta}^2 + \{\frac{g}{l}=9.80\}{\rm{cos}}\theta$ \\
		Stochastic Duffing oscillator & $\mathcal{H}_{s}$=$\frac{1}{2}\dot{X}^2 + \frac{1}{2}\{k=999.35\}X^2 + \frac{1}{4}\{\alpha=2502.31\}X^4$ \\
		Stochastic 3DOF structural system & $\mathcal{H}_{s}$=$\frac{1}{2} \sum_{k=1}^{3} m_k \dot{X}_k^2 + \Bigl( \frac{1}{2}\{k_1=999.99\} X_1^2 + \frac{1}{2}\{k_2=1000.04\}(X_2-X_1)^2 + $ \\
        & $\frac{1}{2}\{k_3=999.96\}(X_3-X_2)^2 \Bigr)$ \\
		Stochastic wave equation & $\mathcal{H}_{s}$=$\frac{1}{2}\sum_{i}{(\partial_t u_i(X,t))^2} + \{\frac{\mu}{\rho}=4.0013\} \sum_{i}{(\partial_x u_i(X,t))^2}$ \\
		Stochastic Euler-Bernoulli Beam & $\mathcal{H}_{s}$=$\frac{1}{2}\sum_{i}{(\partial_t u_i(X,t))^2} + \frac{1}{2} \{\frac{EI}{\rho A}=0.1024\} {(\partial_{xx} u_i(X,t))^2}$ \\
		\bottomrule
	\end{tabular}
\end{table}

\section{Conclusions}\label{sec:conclusion}
The Lagrangian functional stands as a fundamental yet often obscured component in understanding the behavior of dynamical systems, especially those subject to stochastic excitations. In response to this challenge, we propose a novel automated framework for Lagrangian discovery within stochastically excited systems. The proposed framework systematically generates a comprehensive library of terms that may represent the Lagrangian functional.  Leveraging the stochastic Euler-Lagrangian formalism, this approach evaluates and selects appropriate terms while determining their associated coefficients. The proposed approach yields interpretable Lagrangians that can be further utilized for the equation discovery and obtaining the Hamiltonian.
The resulting framework offers a comprehensive exploration of system dynamics under stochastic excitations. It not only provides an interpretable description of the underlying Lagrange density but also identifies the interpretable form of the generalized stochastic force, thereby addressing the limitations of existing approaches to deterministic systems.

The efficacy of our proposed framework is exemplified through six numerical case studies, encompassing stochastic differential equations (SDEs) and stochastic partial differential equations (SPDEs). The outcomes of these tests are compelling and highly encouraging, affirming the robustness of the proposed methodology. In all examples, the discovered Lagrangian density provides an almost exact approximation to the true system.   
The versatility of our framework in handling both discrete (the SDEs) and continuous (the SPDEs) systems with temporal and spatial derivatives broadens the scope of its utility, making it adaptable to a wider array of complex dynamical systems.
Our further investigation into the quality of the derived equations of motion indicates a relative error of less than 1\% in most cases. The predictions using the derived equations of motion almost exactly reproduce the ground truth system responses, indicating the robustness of the proposed framework for the simultaneous discovery of Lagrangian and equations of motion of stochastic dynamical processes. We also investigate the automated discovery of Hamiltonian in stochastic dynamical systems. 
Overall, this is the first time a framework has been proposed for the simultaneous discovery of Lagrangian, Hamiltonian, and equations of motion of stochastic dynamical systems. Our automated framework represents a significant step forward in advancing the understanding and analysis of complex systems continuously subjected to stochastic excitation and whose knowledge is unknown to us a prior.

\section*{Acknowledgements}
T. Tripura acknowledges the financial support received from the Ministry of Education (MoE), India, in the form of the Prime Minister’s Research Fellowship (PMRF). 
B. Hazra gratefully acknowledges the support from SERB, DST India, under Project No. IMP/2019/00276.
S. Chakraborty acknowledges the financial support received from the Science and Engineering Research Board (SERB), India, via grant no. SRG/2021/000467, Ministry of Port and Shipping via letter no. ST14011/74/MT (356529), and seed grant received from IIT Delhi. 

\section*{Code availability}
Upon acceptance, all the data and source codes required to reproduce the results and figures presented in this paper will be made available on GitHub (with a link in the paper).

\section*{Competing interests} 
The authors declare that they have no conflict of interest.

\appendix
\section{Algorithm for sequential threshold least squares}
\begin{algorithm}[H]
    \small
    \caption{Sequential Threshold Least Squares}
    \begin{algorithmic}[1]
        \Require{Mean dictionary matrix $\mathbb{E}[\mathbf{D}_{-u_t^2}]$, regression label $\mathbb{E}[\mathbf{D}_{u_t^2}]$, threshold parameter $\lambda$, maximum number of sequential iterations.}
        \State Obtain initial guess for coefficients ${\bm{C}}_{i} = \arg. \mathop {\min } \limits_{{\bm{C}}^\star_{i}} \left\| \left(\frac{\partial }{\partial u} - \frac{\partial}{\partial t} \frac{\partial }{\partial u_t} \right) \mathbb{E}\left[{\bm{D}}_{-u_t^2}\right]{\bm{C}}_i^{\star \top} - \left(\frac{\partial }{\partial u} - \frac{\partial}{\partial t} \frac{\partial }{\partial u_t} \right) \mathbb{E}\left[{\bm{D}}_{u_t^2}\right]\right\|_2$.
        \For{$k = 1$ \textbf{to} maximum iterations}
            \State{Find the indices $\mathbb{I}(\vert \bm{C}_{ij} \vert < \lambda)$ for $j=1,\ldots,m$.} \label{step:1}
            \State{Set coefficients $\bm{C}_{ij}$ with indices $\mathbb{I}(\vert \bm{C}_{ij} \vert < \lambda)$ equal to zero for $j=1,\ldots,m$.}
            \State{Retain coefficients $\bm{C}_{ij}$ with $\mathbb{I}(\vert \bm{C}_{ij} \vert > \lambda)$ for $j=1,\ldots,m$.} \label{step:2}
            \State Repeat steps \ref{step:1} to \ref{step:2} on the reduced coefficient vector $\bm{C}_i^{\star}$ until indices $\mathbb{I}(\vert \bm{C}_{ij} \vert < \lambda)=0$ for $j=1,\ldots,m$.
            \If{If convergence is occurred,}
                \State{Map the converged coefficients to the original dimension.}
            \EndIf
        \EndFor
        \Ensure{Sparse regression coefficient vector for $\bm{C}_i$.}
    \end{algorithmic}
    \label{algo}
\end{algorithm} 
%


\begin{thebibliography}{10}

\bibitem{mann2018lagrangian}
Peter Mann.
\newblock {\em Lagrangian and Hamiltonian dynamics}.
\newblock Oxford University Press, 2018.

\bibitem{materassi2020stochastic}
Massimo Materassi.
\newblock Stochastic lagrangians for noisy dynamics.
\newblock {\em Chaos, Solitons \& Fractals}, 134:109713, 2020.

\bibitem{finzi2020simplifying}
Marc Finzi, Ke~Alexander Wang, and Andrew~G Wilson.
\newblock Simplifying hamiltonian and lagrangian neural networks via explicit constraints.
\newblock {\em Advances in neural information processing systems}, 33:13880--13889, 2020.

\bibitem{calkin1996lagrangian}
Melvin~G Calkin.
\newblock {\em Lagrangian and Hamiltonian mechanics}.
\newblock World Scientific, 1996.

\bibitem{eyink2020stochastic}
Gregory~L Eyink, Akshat Gupta, and Tamer~A Zaki.
\newblock Stochastic lagrangian dynamics of vorticity. part 1. general theory for viscous, incompressible fluids.
\newblock {\em Journal of Fluid Mechanics}, 901:A2, 2020.

\bibitem{bou2009stochastic}
Nawaf Bou-Rabee and Houman Owhadi.
\newblock Stochastic variational integrators.
\newblock {\em IMA Journal of Numerical Analysis}, 29(2):421--443, 2009.

\bibitem{panda2023geometry}
Satyam Panda, Ankush Gogoi, Budhaditya Hazra, and Vikram Pakrashi.
\newblock Geometry preserving ito-taylor formulation for stochastic hamiltonian dynamics on manifolds.
\newblock {\em Applied Mathematical Modelling}, 119:626--647, 2023.

\bibitem{holm2018stochastic}
Darryl~D Holm and Tomasz~M Tyranowski.
\newblock Stochastic discrete hamiltonian variational integrators.
\newblock {\em BIT Numerical Mathematics}, 58:1009--1048, 2018.

\bibitem{lindgren2019quantum}
Jussi Lindgren and Jukka Liukkonen.
\newblock Quantum mechanics can be understood through stochastic optimization on spacetimes.
\newblock {\em Scientific reports}, 9(1):19984, 2019.

\bibitem{trac2003primer}
Hy~Trac and Ue-Li Pen.
\newblock A primer on eulerian computational fluid dynamics for astrophysics.
\newblock {\em Publications of the Astronomical Society of the Pacific}, 115(805):303, 2003.

\bibitem{li2023lagrangian}
Gang Li and NH~Bian.
\newblock Lagrangian stochastic model for the motions of magnetic footpoints on the solar wind source surface and the path lengths of boundary-driven interplanetary magnetic field lines.
\newblock {\em The Astrophysical Journal}, 945(2):150, 2023.

\bibitem{toth2019hamiltonian}
Peter Toth, Danilo~Jimenez Rezende, Andrew Jaegle, S{\'e}bastien Racani{\`e}re, Aleksandar Botev, and Irina Higgins.
\newblock Hamiltonian generative networks.
\newblock {\em arXiv preprint arXiv:1909.13789}, 2019.

\bibitem{greydanus2019hamiltonian}
Samuel Greydanus, Misko Dzamba, and Jason Yosinski.
\newblock Hamiltonian neural networks.
\newblock {\em Advances in neural information processing systems}, 32, 2019.

\bibitem{sanchez2019hamiltonian}
Alvaro Sanchez-Gonzalez, Victor Bapst, Kyle Cranmer, and Peter Battaglia.
\newblock Hamiltonian graph networks with ode integrators.
\newblock {\em arXiv preprint arXiv:1909.12790}, 2019.

\bibitem{lutter2019deep}
Michael Lutter, Kim Listmann, and Jan Peters.
\newblock Deep lagrangian networks for end-to-end learning of energy-based control for under-actuated systems.
\newblock In {\em 2019 IEEE/RSJ International Conference on Intelligent Robots and Systems (IROS)}, pages 7718--7725. IEEE, 2019.

\bibitem{cranmer2020lagrangian}
Miles Cranmer, Sam Greydanus, Stephan Hoyer, Peter Battaglia, David Spergel, and Shirley Ho.
\newblock Lagrangian neural networks.
\newblock {\em arXiv preprint arXiv:2003.04630}, 2020.

\bibitem{gruver2022deconstructing}
Nate Gruver, Marc Finzi, Samuel Stanton, and Andrew~Gordon Wilson.
\newblock Deconstructing the inductive biases of hamiltonian neural networks.
\newblock {\em arXiv preprint arXiv:2202.04836}, 2022.

\bibitem{bhattoo2022learning}
Ravinder Bhattoo, Sayan Ranu, and NM~Krishnan.
\newblock Learning articulated rigid body dynamics with lagrangian graph neural network.
\newblock {\em Advances in Neural Information Processing Systems}, 35:29789--29800, 2022.

\bibitem{tripura2023bayesian}
Tapas Tripura and Souvik Chakraborty.
\newblock A bayesian framework for discovering interpretable lagrangian of dynamical systems from data.
\newblock {\em arXiv preprint arXiv:2310.06241}, 2023.

\bibitem{tripura2024discovering}
Tapas Tripura and Souvik Chakraborty.
\newblock Discovering interpretable lagrangian of dynamical systems from data.
\newblock {\em Computer Physics Communications}, 294:108960, 2024.

\bibitem{schmidt2009distilling}
Michael Schmidt and Hod Lipson.
\newblock Distilling free-form natural laws from experimental data.
\newblock {\em science}, 324(5923):81--85, 2009.

\bibitem{brunton2016discovering}
Steven~L Brunton, Joshua~L Proctor, and J~Nathan Kutz.
\newblock Discovering governing equations from data by sparse identification of nonlinear dynamical systems.
\newblock {\em Proceedings of the national academy of sciences}, 113(15):3932--3937, 2016.

\bibitem{nayek2021spike}
Rajdip Nayek, Ramon Fuentes, Keith Worden, and Elizabeth~J Cross.
\newblock On spike-and-slab priors for bayesian equation discovery of nonlinear dynamical systems via sparse linear regression.
\newblock {\em Mechanical Systems and Signal Processing}, 161:107986, 2021.

\bibitem{nayek2022equation}
Rajdip Nayek, Keith Worden, and Elizabeth~J Cross.
\newblock Equation discovery using an efficient variational bayesian approach with spike-and-slab priors.
\newblock In {\em Model Validation and Uncertainty Quantification, Volume 3: Proceedings of the 39th IMAC, A Conference and Exposition on Structural Dynamics 2021}, pages 149--161. Springer, 2022.

\bibitem{wentz2023derivative}
Jacqueline Wentz and Alireza Doostan.
\newblock Derivative-based sindy (dsindy): Addressing the challenge of discovering governing equations from noisy data.
\newblock {\em Computer Methods in Applied Mechanics and Engineering}, 413:116096, 2023.

\bibitem{schaeffer2017learning}
Hayden Schaeffer.
\newblock Learning partial differential equations via data discovery and sparse optimization.
\newblock {\em Proceedings of the Royal Society A: Mathematical, Physical and Engineering Sciences}, 473(2197):20160446, 2017.

\bibitem{rudy2017data}
Samuel~H Rudy, Steven~L Brunton, Joshua~L Proctor, and J~Nathan Kutz.
\newblock Data-driven discovery of partial differential equations.
\newblock {\em Science advances}, 3(4):e1602614, 2017.

\bibitem{chen2021physics}
Zhao Chen, Yang Liu, and Hao Sun.
\newblock Physics-informed learning of governing equations from scarce data.
\newblock {\em Nature communications}, 12(1):6136, 2021.

\bibitem{flaschel2021unsupervised}
Moritz Flaschel, Siddhant Kumar, and Laura De~Lorenzis.
\newblock Unsupervised discovery of interpretable hyperelastic constitutive laws.
\newblock {\em Computer Methods in Applied Mechanics and Engineering}, 381:113852, 2021.

\bibitem{joshi2022bayesian}
Akshay Joshi, Prakash Thakolkaran, Yiwen Zheng, Maxime Escande, Moritz Flaschel, Laura De~Lorenzis, and Siddhant Kumar.
\newblock Bayesian-euclid: Discovering hyperelastic material laws with uncertainties.
\newblock {\em Computer Methods in Applied Mechanics and Engineering}, 398:115225, 2022.

\bibitem{more2023bayesian}
Kalpesh~Sanjay More, Tapas Tripura, Rajdip Nayek, and Souvik Chakraborty.
\newblock A bayesian framework for learning governing partial differential equation from data.
\newblock {\em Physica D: Nonlinear Phenomena}, 456:133927, 2023.

\bibitem{boninsegna2018sparse}
Lorenzo Boninsegna, Feliks N{\"u}ske, and Cecilia Clementi.
\newblock Sparse learning of stochastic dynamical equations.
\newblock {\em The Journal of chemical physics}, 148(24), 2018.

\bibitem{li2021data}
Yang Li and Jinqiao Duan.
\newblock A data-driven approach for discovering stochastic dynamical systems with non-gaussian l{\'e}vy noise.
\newblock {\em Physica D: Nonlinear Phenomena}, 417:132830, 2021.

\bibitem{tripura2023sparse}
Tapas Tripura and Souvik Chakraborty.
\newblock A sparse bayesian framework for discovering interpretable nonlinear stochastic dynamical systems with gaussian white noise.
\newblock {\em Mechanical Systems and Signal Processing}, 187:109939, 2023.

\bibitem{tripura2023probabilistic}
Tapas Tripura, Aarya~Sheetal Desai, Sondipon Adhikari, and Souvik Chakraborty.
\newblock Probabilistic machine learning based predictive and interpretable digital twin for dynamical systems.
\newblock {\em Computers \& Structures}, 281:107008, 2023.

\bibitem{mathpati2024discovering}
Yogesh~Chandrakant Mathpati, Tapas Tripura, Rajdip Nayek, and Souvik Chakraborty.
\newblock Discovering stochastic partial differential equations from limited data using variational bayes inference.
\newblock {\em Computer Methods in Applied Mechanics and Engineering}, 418:116512, 2024.

\bibitem{Noether1}
Emmy Noether.
\newblock Invariant variation problems.
\newblock {\em Transport Theory and Statistical Physics}, 1(3):186--207, 1971.

\bibitem{kaheman2020sindy}
Kadierdan Kaheman, J~Nathan Kutz, and Steven~L Brunton.
\newblock Sindy-pi: a robust algorithm for parallel implicit sparse identification of nonlinear dynamics.
\newblock {\em Proceedings of the Royal Society A}, 476(2242):20200279, 2020.

\bibitem{ito1944109}
Kiyosi It{\^o}.
\newblock 109. stochastic integral.
\newblock {\em Proceedings of the Imperial Academy}, 20(8):519--524, 1944.

\bibitem{tankov2003financial}
Peter Tankov.
\newblock {\em Financial modelling with jump processes}.
\newblock CRC press, 2003.

\bibitem{turelli1977random}
Michael Turelli.
\newblock Random environments and stochastic calculus.
\newblock {\em Theoretical population biology}, 12(2):140--178, 1977.

\bibitem{chechkin2014marcus}
Alexei Chechkin and Ilya Pavlyukevich.
\newblock Marcus versus stratonovich for systems with jump noise.
\newblock {\em Journal of Physics A: Mathematical and Theoretical}, 47(34):342001, 2014.

\bibitem{williams2012stochastic}
Brett~T Williams.
\newblock {\em On stochastic differential equations in the Ito and in the Stratonovich sense}.
\newblock University of Maryland, College Park, 2012.

\bibitem{kloeden1992higher}
Peter~E Kloeden and Eckhard Platen.
\newblock Higher-order implicit strong numerical schemes for stochastic differential equations.
\newblock {\em Journal of statistical physics}, 66(1-2):283--314, 1992.

\bibitem{duong2021hamiltonian}
Thai Duong and Nikolay Atanasov.
\newblock Hamiltonian-based neural ode networks on the se (3) manifold for dynamics learning and control.
\newblock {\em arXiv preprint arXiv:2106.12782}, 2021.

\end{thebibliography}

\newcommand{\noopsort}[1]{} \newcommand{\printfirst}[2]{#1} \newcommand{\singleletter}[1]{#1} \newcommand{\switchargs}[2]{#2#1}

\end{document}